\documentclass[12pt]{amsart}
\usepackage{amsmath,amsfonts,euscript,amscd,amsthm,amssymb,upref,graphics}

\theoremstyle{plain}
\newtheorem{theorem}{Theorem}
\swapnumbers

\newtheorem{proposition}[subsection]{Proposition}
\newtheorem{lemma}[subsection]{Lemma}
\newtheorem{corollary}[subsection]{Corollary}

\theoremstyle{definition}
\newtheorem{definition}[subsection]{Definition}
\newtheorem{example}[subsection]{Example}
\newtheorem{remark}[subsection]{Remark}

\newtheorem{nothing*}[subsection]{}
\newtheorem{notation}[subsection]{Notation}
\newtheorem{convention}[subsection]{Convention}

\newcommand{\rien}[1]{}
\newcommand{\Span}{ \operatorname{{\rm Span}}}

\newcommand{\Image}{ \operatorname{{\rm Im}}}
\newcommand{\diver}{ \operatorname{{\rm div}}}
\newcommand{\Lie}{ \operatorname{{\rm Lie}}}
\newcommand{\AVF}{ \operatorname{{\rm AVF}}}
\newcommand{\VFH}{ \operatorname{{\rm VF}_{hol}}}

\newcommand{\LieH}{ \operatorname{{\rm Lie}_{hol}}}
\newcommand{\VFA}{ \operatorname{{\rm VF}_{alg}}}
\newcommand{\VF}{ \operatorname{{\rm VF}}}
\newcommand{\LieA}{ \operatorname{{\rm Lie}_{alg}}}

\newcommand{\IVFA}{ \operatorname{{\rm IVF}_{alg}}}

\newcommand{\VFHO}{ \operatorname{{\rm VF}_{hol}^\omega}}
\newcommand{\LieHO}{ \operatorname{{\rm Lie}_{hol}^\omega}}
\newcommand{\VFAO}{ \operatorname{{\rm VF}_{alg}^\omega}}
\newcommand{ \LieAO}{ \operatorname{{\rm Lie}_{alg}^\omega}}

\newcommand{\IVFAO}{ \operatorname{{\rm IVF}_{alg}^\omega}}

\newcommand{\IVF}{ \operatorname{{\rm IVF}}}
\newcommand{\Div}{ \operatorname{{\rm Div}}}

\newcommand{\pr}{ \operatorname{{\rm pr}}}
\newcommand{\av}{ \operatorname{{\rm av}}}

\newcommand{\C}{\ensuremath{\mathbb{C}}}
\newcommand{\Z}{\ensuremath{\mathbb{Z}}}

\newcommand{\R}{\ensuremath{\mathbb{R}}}
\newcommand{\T}{\ensuremath{\mathbb{T}}}

\newcommand{\sgoth}{{\ensuremath{\mathfrak{s}}}}
\newcommand{\lgoth}{{\ensuremath{\mathfrak{l}}}}

\newcommand{\ggoth}{{\ensuremath{\mathfrak{g}}}}
\newcommand{\hgoth}{{\ensuremath{\mathfrak{h}}}}

\newcommand{\cL}{{\ensuremath{\mathcal{L}}}}

\newcommand{\cF}{{\ensuremath{\mathcal{F}}}}
\newcommand{\cS}{{\ensuremath{\mathcal{S}}}}

\newcommand{\cT}{{\ensuremath{\mathcal{T}}}}
\newcommand{\cZ}{{\ensuremath{\mathcal{Z}}}}

\newcommand{\Ker}{{\rm Ker} \,}

\def\vfo{{\rm VF}_{hol}}

\renewcommand{\epsilon}{\varepsilon}
\renewcommand{\phi}{\varphi}

\begin{document}
\renewcommand{\baselinestretch}{1.07}

\title[Algebraic Volume Density Property of Affine Algebraic Manifolds]
{Algebraic Volume Density Property of Affine Algebraic Manifolds}
\author{Shulim Kaliman}
\address{Department of Mathematics\\
University of Miami\\
Coral Gables, FL 33124 \ \ USA}
\email{kaliman@math.miami.edu}
\author{Frank Kutzschebauch}
\address{Mathematisches Institut \\Universit\"at Bern
 \\Sidlerstr. 5
 \\ CH-3012 Bern, Switzerland}
\thanks{{\bf Acknowledgements:} This research was started during a visit
of the first author to the University of Bern  and continued
during a visit of the second author to the University of Miami,
Coral Gables. We thank these institutions for their generous
support and excellent working conditions. The research of the
second author was also partially supported by Schweizerische
Nationalfonds grant No 200020-124668 / 1 .}
\email{Frank.Kutzschebauch@math.unibe.ch} 
{\renewcommand{\thefootnote}{} }
\begin{abstract} We introduce the notion of algebraic volume density property for affine algebraic manifolds and prove some important basic facts about it, in particular that it implies the volume density property. The main results of the paper are producing two big classes of examples of Stein manifolds with volume density property. One class consists of certain affine modifications of $\C^n$ equipped with a canonical volume form, the other is the class of all Linear Algebraic Groups equipped with the left invariant volume form.
\end{abstract}
\maketitle \vfuzz=2pt

\vfuzz=2pt
\section{Introduction}

In this paper we study a less developed part of  the Anders\' en-Lempert
theory  (\cite{A}, \cite{AL}, \cite{FR}, \cite{V1}, \cite {TV1}, \cite{TV2}, \cite{F}) namely the case of volume preserving maps. Recall that Anders\' en-Lempert theory describes complex manifolds such that among other things  the local phase flows on their holomorphically convex compact
subsets can be approximated by global holomorphic
automorphisms which leads to construction of  holomorphic
automorphisms with prescribed local properties. Needless to say that
this implies some remarkable consequences for such manifolds
(e.g., see \cite{V1}, \cite{V2}, \cite{KaKu1}). It turns out that a complex
manifold has such approximations if it possesses the following density property introduced by \textsc{Varolin}.

\begin{definition}\label{int.1}
A complex manifold $X$ has the density property if in the
compact-open topology the Lie algebra $\LieH (X)$ generated by
completely integrable holomorphic vector fields on $X$ is dense in
the Lie algebra $\VFH (X)$ of all holomorphic vector fields on
$X$. An affine algebraic manifold $X$ has the algebraic density
property if the Lie algebra $\LieA (X)$ generated by completely
integrable algebraic vector fields on it coincides with  the Lie
algebra $\VFA (X)$ of all algebraic vector fields on it (clearly,
the algebraic density property implies the density property).
\end{definition}
The algebraic density property was established for a wide variety
of affine algebraic manifolds, including all connected linear
algebraic groups except for $\C_+$ and complex tori by the authors \cite{KaKu1}.
Furthermore, in the coming paper of \textsc{Donzelli, Dvorsky} and the first author \cite{DoDvKa} it will be extended to homogeneous affine algebraic manifolds
different from $\C_+$, complex tori, and one extra surface.

However \textsc{Anders\'en, Lempert, Forstneric, Rosay} and \textsc{Varolin} considered also another property which has similar consequences for automorphisms preserving a
volume form.
\begin{definition}\label{int.2}
Let a complex manifold  $X$ be
equipped with a holomorphic volume form $\omega$
(i.e. $\omega$ is nowhere vanishing section of the canonical
bundle). We say that $X$ has the volume density property with
respect to $\omega$ if in the compact-open topology the Lie
algebra $\LieHO$ generated by globally integrable holomorphic
vector fields $\nu$ such that $\nu (\omega)=0$, is dense in the
Lie algebra $\VFHO (X)$ of all holomorphic vector fields that
annihilate $\omega$ (note that condition $\nu (\omega)=0$ is
equivalent to the fact that $\nu$ is of zero $\omega$-divergence).

\end{definition}

Compared with the density property, the class
of complex manifolds with established volume density property has been
quite narrow. It was
essentially described by the original result of
\textsc{Anders\'en} and \textsc{Lempert} \cite{A}, \cite{AL} who proved it for
Euclidean spaces plus a few other examples found by
\textsc{Varolin} \cite{V1}. In particular he proved that $SL_2 (\C)$ has volume density property with respect
to the Haar form but he was unable to decide whether the following hypersurface given by a similar
equation like $SL_2 (\C)$

$$ \Sigma^3 = \{ (a,b,c,d) \in \C^4 \ : \ a^2 c - b d = 1\}$$
had volume density property or not (\cite{V3} section 7).

In order to deal with this lack of examples we introduce
like in  the previous pattern the following.

\begin{definition}\label{int.3}
If $X$ is affine algebraic we say that $X$ has the algebraic
volume density property with respect to an algebraic volume form $\omega$ if the Lie
algebra $\LieAO$ generated by globally integrable algebraic vector
fields $\nu$ such that $\nu (\omega)=0$, coincides with the Lie
algebra $\VFAO (X)$ of all algebraic vector fields that annihilate
$\omega$.
\end{definition}

It is much more difficult to establish  the algebraic volume density property than
the algebraic density property. This is caused, perhaps, by
the following difference which does not allow to apply
the most effective criterion for the algebraic density property (see \cite{KaKu2}):
$\VFAO (X)$ is not a module over
the ring $\C [X ]$ of regular functions on $X$ while $\VFA (X)$ is.
Furthermore, some features that are straightforward for the algebraic density property
are not at all clear in the volume-preserving case. For instance, it
is not quite obvious that the algebraic volume density property
implies the volume density property and that the product of two manifolds with algebraic volume
density property has again the algebraic volume density property.
We shall show in this paper the validity of these two facts
among other results  that enable us to
enlarge the class of examples of Stein manifolds with the volume density property substantially.
In particular we establish the following.

\begin{theorem}\label{one} Let $X'$ be a hypersurface in $\C_{u,v,\bar x}^{n+2}$ given by an equation
of form $P(u,v,\bar x )=uv-p(\bar x )=0$ where $p$ is a polynomial on $\C^n_{\bar x}$ with a smooth reduced zero fiber $C$ such that  its reduced cohomology $\hat H^{n-2}(C, \C )=0$.  Let $\Omega$ be a standard volume form on $\C^{n+2}$ and $\omega'$ be a volume form on $X'$ such that $\omega' \wedge dP =\Omega |_{X'}$.  Then $X'$ has the algebraic $\omega'$-volume density property.

\end{theorem}

This gives, of course, an affirmative answer to \textsc{Varolin}'s question mentioned before. The next
theorem is
our main result.

\begin{theorem}\label{theorem} Let $G$ be a linear algebraic
 group. Then $G$ has the algebraic
volume density property with respect to the left (or right) invariant volume
form.
\end{theorem}

Let us describe briefly the content of the paper and the main steps in the proof of these facts.

In Section 2
we remind some standard facts about divergence.

In Section 3
we deal with Theorem \ref{one} in a slightly more general situation. Namely we consider a hypersurface
$X'$ in $X\times \C^2_{u,v}$ given by an equation $P:= uv-p(x)=0$ where $X$ is a smooth affine
algebraic variety and $p(x)$ is a regular function on $X$. We suppose that $X$ is equipped with
a volume form $\omega$ and establish the existence of a volume form $\omega'$ on $X'$ such that
$\Omega |_{X'} =dP \wedge \omega$ where $\Omega = du \wedge dv \wedge \omega$. Then we prove
(Proposition \ref{A}) that $X'$ has the $\omega'$-volume algebraic density property provided
two technical conditions (A1) and (A2) hold.

Condition (A2) is easily verifiable for $X'$ and condition (A1) is equivalent to the following (Lemma \ref{A'}):
the space of zero $\omega$-divergence algebraic vector fields on $X$
tangent to the zero fiber $C$ of $p$ is generated by vector fields of form $\nu_1 (fp) \nu_2-
\nu_2 (fp) \nu_1$ where $\nu_1$ and $\nu_2$ are commutative completely integrable algebraic
vector fields of zero $\omega$-divergence on $X$ and $f$ is a regular function on $X$.

Then we notice the duality between the spaces of zero $\omega$-divergence vector fields on $X$ and
closed $(n-1)$-forms on $X$ which is achieved via the inner product that assigns
to each vector field $\nu$ the $(n-1)$-form $\iota_\nu (\omega )$ (Lemma \ref{add.5}).
This duality allows to reformulate condition (A1) as the following:

(i) the space of
algebraic $(n-2)$-forms on $X$ is generated by the forms of type $\iota_{\nu_1} \iota_{\nu_2} (\omega )$
where $\nu_1$ and $\nu_2$ are as before; and

(ii) the outer differentiation sends the space of $(n-2)$-forms on $X$ that vanish on $C$ to the set
of $(n-1)$-forms whose restriction to $C$ yield the zero $(n-1)$-form on $C$.

In the case
of $X$ isomorphic to a Euclidean space (i) holds automatically with $\nu_1$ and $\nu_2$ running
over the set of partial derivatives.

If the reduced cohomology $\hat H^{n-2} (C, \C )=0$
and also $H^n(X, \C )=0$ the
validity of (ii) is a consequence of the Grothendieck theorem (see Proposition \ref{sur.1})
that states that the complex
cohomology can be computed via the De Rham complex of algebraic forms on a smooth affine
algebraic variety which concludes the proof of Theorem \ref{one}.

We end Section 3 with
an important corollary of Theorem \ref{one} which will be important for the proof of Theorem \ref{theorem} : the groups $SL_2(\C )$ (already proved by \textsc{Varolin} as mentioned above) and $PSL_2(\C )$ have the algebraic volume density property with respect to the invariant volume  (Propositions \ref{sl2'} and \ref{psl2}). The proof is based on the fact that $SL_2(\C )$ is isomorphic to
the hypersurface in $\C^4_{u,v,x_1,x_2}$ given by $uv-x_1x_2-1=0$.

Section 4 contains two general facts about the algebraic volume density property with short but
non-trivial proofs. The first of them (Proposition \ref{approx})
says that the algebraic volume density property
implies the volume density property (in the holomorphic sense). It is also based on the Grothendieck
theorem mentioned before. The second one (Proposition \ref{pro1})
states that the product $X\times Y$ of two affine algebraic
manifolds $X$ and $Y$ with the algebraic volume density property
(with respect to volumes $\omega_X$ and $\omega_Y$) has also the algebraic volume density
property (with respect to $\omega_X \times \omega_Y$). As a consequence of this result we establish
the algebraic volume density property for all tori which was also established
earlier by \textsc{Varolin} \cite{V1} (recall that the density property is
not established for higher dimensional tori yet and the algebraic density property does not
hold for these objects \cite{A'}).

We start Section 5 
with discussion of a phenomenon
which makes the proof of Proposition \ref{pro1} about the algebraic volume density of  $X \times Y$
non-trivial and prevents us from spreading it directly to locally trivial fibrations.
More precisely, consider the subspace $F_Y$ of $\C [Y]$ generated by the images
${\rm Im} \, \delta$ with $\delta$
running over $\LieAO (Y)$.   In general $F_Y \ne \C [Y]$ and the absence of equality here
is the source of difficulties.
Nevertherless
one can follow the pattern of the proof of Proposition \ref{pro1} in the
case of fibrations when the span of $F_Y$ and constants yields $\C [Y]$.
A manifold $Y$ with this property is called fine and
we describe simple facts about such objects which include $SL_2$ and $PSL_2$.
We also introduce the notion of a refined
volume fibration $p : W \to X$ which a generalization of the product situation.
Among other assumptions the fiber  of $p$, the base $X$, and the total space $W$ are equipped with
nicely related volume forms and the fiber is a fine manifold with the algebraic density
property.
The main result in Section 5 is Theorem  \ref{ffi.th1}
saying that the total space of
a refined volume fibration  has the algebraic volume density property
provided  its  base has the algebraic volume density property as well.

Section 6 contains basic knowledge about invariant volume forms on linear algebraic groups.
Of further  importance will be Corollary \ref{refine.2} about the Mostow decomposition of a linear algebraic group as the product of Levi reductive subgroup and it's unipotent radical. We end it with
an important example of a refined volume fibration - the quotient map of a reductive group by its Levi semi-simple subgroup (see Lemma \ref{refine.1}).

Section 7 prepares the proof of Theorem \ref{theorem} in the case of a semi-simple group. The central
notion discussed in that section is a $p$-compatible vector field  $\sigma' \in \LieAO (W)$ for a locally trivial fibration
$p : W \to X$. Its most important property is that ${\rm Span} \, \Ker \sigma' \cdot \Ker \delta'$ coincides
with the algebra $\C [W]$ of regular functions for any $\delta' \in \VFAO (W)$ tangent to
the fibers of $p$. It is established that for any at least three-dimensional semi-simple group $G$
and its $SL_2$- or $PSL_2$-subgroup $S$ corresponding to the root of the Dynkin diagram
the fibration $q : G \to G/S$ admits a sufficiently large family $q$-compatible vector fields.
The existence of such a family is an assumption in the definition of a refined fibration which
(in combination with the fact that $SL_2$ and $PSL_2$ are fine and have
 the algebraic volume density property)
leads to the claim that $q$ is refined. This enables us to use properties of refined fibrations established
earlier in Proposition \ref{vfi.pro2} (but not  Theorem  \ref{ffi.th1} since it is unknown
whether $G/S$ has the algebraic volume density property).

Section 8 contains the proof of Theorem \ref{theorem}. The general case follows easily
from a semi-simple one (via  Lemma \ref{refine.1}, Theorem \ref{ffi.th1}, and Corollary \ref{refine.2}). The idea of the proof in the latter case is the following.
We consider  $SL_2$- or $PSL_2$-subgroups $S_0, \ldots , S_m$ corresponding
to the simple roots of a Dynkin diagram of a semi-simple group $G$ and
fibrations $p_i : G\to G/S_i$ with $ i =0, \ldots , m$.  Using results of Section 7 we
establish that there is  a sufficiently big collection $\Theta$
of completely integrable fields $\theta $ of zero divergence that are of $p_i$-compatible
for every $i$.  Furthermore up to an element of $\LieAO (G)$ every algebraic vector field of zero divergence can
be presented as a finite sum $\sum h_i \theta_i$ where $\theta_i \in \Theta$ and
$h_i \in \C [G]$. Then we consider a standard averaging operator $\av_j$ on $\C [G]$ that
assigns to each $h \in \C [G]$ a regular function $\av_j (h)$ invariant with respect
to the natural $S_j$-action on $G$ and establish the following relation:
 $\sum h_i \theta_i\in \LieAO (G)$  if and only if  $\sum \av_j ( h_i) \theta_i \in \LieAO (G)$.
 We show also that a consequent application of operators $\av_0, \ldots , \av_m$ leads to
 a function invariant with respect to each $S_j , \, j=0, \ldots , m$. Since the only
 functions invariant under the natural actions of all such subgroups are constants
 we see that   $\sum h_i \theta_i\in \LieAO (G)$  because $\sum c_i \theta_i \in \LieAO (G)$
 for constant coefficients $c_i$ which concludes the proof of Theorem \ref{theorem}.

The appendix contains definition of strictly semi-compatible fields and  refinements of two Lemmas about it from
our previous work \cite{KaKu1}.

{\em Acknowledgments.} We would like to thank \textsc{A. Dvorsky} for
helpful consultations.

\section{Preliminaries}\label{Preliminaries}

Recall that a holomorphic vector field $V\in \vfo (\C^n)$ is
completely (or globally) integrable if for any initial value $z
\in \C^n$ there is a global holomorphic solution of the
ordinary differential equation
\begin{equation}\label{flow}
\dot \gamma (t)= V(\gamma (t)), \quad \gamma(0)=z.
\end{equation}
In this case the phase flow (i.e. the map $\C\times \C^n \to
\C^n$ given by $(t,z)\mapsto \gamma_z (t)$) is a holomorphic
action of the additive group $\C_+$ on $\C^n$, where index
$z$ in $\gamma_z$ denotes the dependence on the initial value. It
is worth mentioning that this action is not necessarily algebraic
in the case of an algebraic vector field $V \in \VFA (\C^n)$.

A holomorphic (resp. algebraic) volume form on
a complex (resp. affine algebraic) manifold $X$ of dimension $n$ is
a nowhere vanishing holomorphic (resp. algebraic) $n$-form.
Let us discuss some simple properties of the divergence
$\diver_{\omega} (\nu )$ of a vector field $\nu$ on $X$ with
respect to this volume form $\omega$. The divergence is defined by
the equation

\begin{equation}\label{DefDivergence}
\diver_{\omega} (\nu ) \omega =L_{\nu} (\omega )
\end{equation}

\noindent where $L_{\nu}$ is the Lie derivative. Furthermore, for any vector
fields $\nu_1,\nu_2$ on $X$ we have the following relation between
divergence and Lie bracket

\begin{equation}\label{BracketDivergence}
\diver_{\omega} ([\nu_1,\nu_2])=L_{\nu_1} (\diver_{\omega}
(\nu_2)) - L_{\nu_2}( \diver_{\omega} (\nu_1)).
\end{equation}

In particular, when $\diver_{\omega} (\nu_1)=0$ we have

\begin{equation}\label{ZeroDivergence}
\diver_{\omega} ([\nu_1,\nu_2])=L_{\nu_1} (\diver_{\omega}
(\nu_2)).
\end{equation}

Another useful formula is

\begin{equation}\label{FactorDivergence}
\diver_{\omega} (f\nu )=f\diver_{\omega} (\nu) + \nu (f)
\end{equation}

for any holomorphic function $f$ on $X$.

\begin{lemma}\label{Divergence}
Let $Y$ be a Stein complex manifold with a volume form $\Omega$ on
it, and $X$ be a submanifold of $Y$ which is a strict complete
intersection (that is, the defining ideal of $X$ is generated by
holomorphic functions $P_1, \ldots , P_k$ on $Y$, where $k$ is the
codimension of $X$ in $Y$). Suppose that $\nu$ is a vector field
on $X$ and $\mu$ is its extension to $Y$ such that $\mu (P_i)=0$
for every $i=1, \ldots , k$. Then

{\rm (i)} there exists a volume form $\omega$ on $X$ such that
$\Omega |_X = dP_1 \wedge \ldots \wedge dP_k \wedge \omega$; and

{\rm (ii)} $\diver_{\omega} (\nu )= \diver_{\Omega } (\mu )|_{X}$.

\end{lemma}

\begin{proof}
Let $x_1, \ldots , x_n$ be a local holomorphic coordinate system
in a neighborhood of a point in $X$. Then $P_1, \ldots , P_k, x_1,
\ldots , x_n$ is a local holomorphic coordinate system in a
neighborhood of this point in $Y$. Hence in that neighborhood
$\Omega = h dP_1\wedge \ldots \wedge dP_k \wedge dx_1 \wedge
\ldots \wedge dx_n$ where $h$ is a holomorphic function. Set
$\omega = h|_X dx_1 \wedge \ldots \wedge dx_n$. This is the
desired volume form in (i).

Recall that $L_{\nu }=d \circ \imath_{\nu} + \imath_{\nu }\circ d$
where $\imath_{\nu}$ is the inner product with respect to $\nu$
(\cite{KoNo}, Chapter 1, Proposition 3.10). Since $\mu (P_i)=0$ we
have  $L_{\mu} (dP_i) =0$. Hence by formula (\ref{DefDivergence})
we have
$\diver_{\Omega}( \mu ) \Omega |_X =$
$$ L_{\mu } \Omega |_X
=L_{\mu} (dP_1 \wedge \ldots \wedge dP_k \wedge \omega )|_X = dP_1
\wedge \ldots \wedge dP_k|_X \wedge L_{\nu} \omega + L_{\mu}(dP_1
\wedge \ldots \wedge dP_k )|_X \wedge \omega $$  $$=  \diver_{\omega}
(\nu) (dP_1 \wedge \ldots \wedge dP_k)|_X \wedge \omega
=\diver_{\omega} (\nu ) \Omega |_X$$ which is (ii).

\end{proof}

\begin{remark}\label{add.1} \hfill

(1) Lemma \ref{Divergence} remains valid
in the algebraic category

(2) Furthermore, it enables us  to
compute the divergence of a vector field on $X$ via the divergence
of a vector field extension on an ambient space. It is worth
mentioning that there is another simple way to compute divergence
on $X$ which leads to the same formulas in Lemma \ref{Facts}
below. Namely, $X$ that we are going to consider will be an affine
modification $\sigma : X \to Z$ of another affine algebraic
manifold $Z$ with a volume form $\omega_0$ (for definitions of
affine and pseudo-affine modifications see \cite{KZ} ). In particular, for some divisors $D \subset Z$
and $E \subset X$ the restriction of $\sigma$ produces an
isomorphism $X \setminus E \to Z \setminus D$. One can suppose
that $D$ coincides with the zero locus of a regular (or
holomorphic) function $\alpha$ on $Z$. In the situation, we are
going to study, the function ${\tilde \alpha}= \alpha \circ
\sigma$ has simple zeros on $E$. Consider the form $\sigma^*
\omega_0$ on $X$. It may vanish on $E$ only. Dividing this form by
some power ${\tilde \alpha}^k$ we get a volume form on $X$. In
order to compute divergence of a vector field on $X$ it suffices
to find this divergence on the Zariski open subset $X \setminus E
\simeq Y \setminus D$, i.e. we need to compute the divergence of a
vector field $\nu$ on $Y\setminus D$ with respect to a volume form
$\beta \omega_0$ where $\beta = \alpha^{-k}$. The following
formula relates it with the divergence with respect to $\omega_0$:
\begin{equation}\label{add.3} \diver_{\beta \omega_0} (\nu ) =\diver_{\omega_0} (\nu) +
L_{\nu } (\beta)/ \beta .\end{equation}

In the cases, we need to consider, $\beta$ will be often in the
kernel of $\nu$, i.e. $\diver_{\beta \omega_0} (\nu )
=\diver_{\omega_0} (\nu)$ in these cases.

(3) If the normal bundle of $X\subset \C^n$ is trivial we may
choose $\omega$ as the restriction of the standard volume form on
$\C^n$ by Lemma \ref{Divergence}.  Indeed, taking $n$
sufficiently large we can always assume that $X$ is a complete
intersection in $\C^n$ (see for example \cite{Sch}).
\end{remark}

The condition in Lemma \ref{Divergence} that an algebraic field $\nu$
on $X$ has an extension $\mu$ on $Y$ with $\mu (P_i)=0$ is also
very mild. We consider it in the case of hypersurfaces only.

\begin{lemma}\label{extension} Let $X$ be a smooth hypersurface in
a complex Stein (resp. affine algebraic) manifold $Y$ given by
zero of a reduced holomorphic (resp. algebraic) function $P$ on
$Y$. Then every holomorphic (resp. algebraic) vector field $\nu$
on $X$ has a similar extension $\mu$ to $Y$ such that $\mu (P)=0$.

\end{lemma}

\begin{proof} Consider, for instance the algebraic case, i.e.
$P$ belongs to the ring $\C [Y]$ of regular functions on $Y$.
Since $\mu$ must be tangent to $X$ we see that $\mu (P)$ vanishes
on $X$, i.e. $\mu (P)=PQ$ where $Q \in \C [Y]$. Any other
algebraic extension of $\nu$ is of form $\tau = \mu - P \theta$
where $\theta \in \VFA (Y)$. Thus if $\theta (P)= Q$ then we are
done.

In order to show that such $\theta$ can be found consider the set
$M =\{ \theta (P) | \theta \in \VFA (Y) \}$. One can see that $M$
is an ideal of $\C [Y]$. Therefore, it generates a coherent
sheaf $\cF$ over $Y$. The restriction $Q|_{Y \setminus X}$ is a
section of $\cF |_{Y \setminus X}$ because $Q= \mu (P)/P$. Since
$X$ is smooth for every point $x \in X$ there are a Zariski open
neighborhood $U$ in $Y$ and  an algebraic vector field $\partial$ such
that $\partial (P)$ does not vanish on $U$. Hence $Q|_U$ is a
section of $\cF |_U$. Since $\cF$ is coherent this implies that
$Q$ is a global section of $\cF$ and, therefore, $Q \in M$ which
is the desired conclusion.
\end{proof}

\begin{nothing*}\label{Terminology} {\bf Terminology and Notation.}
{\rm In the rest of this section $X$ is a closed affine algebraic
submanifold of $\C^n$, $\omega$ is an algebraic volume form on $X$, $p$ is
a regular function on $X$ such that the divisor $p^*(0)$ is smooth
reduced, $X'$ is the hypersurface in $Y= \C^2_{u,v} \times X$
given by the equation $P:=uv -p=0$.\footnote{By abusing notation
we treat $p$ in this formula as a function on $Y$, and, if
necessary, we treat it as a function on $X'$. Furthermore, by
abusing notation, for any regular function on $X$ we denote its
lift-up to $Y$ or $X'$ by the same symbol.} Note that $X'$ is
smooth and, therefore, Lemma \ref{extension} is applicable. We
shall often use the fact that every regular function $f$ on $X'$
can be presented uniquely as the restriction of a regular function
on $Y$ of the form

\begin{equation}\label{StandardForm}
f=\sum_{i=1}^m (a_iu^i+b_iv^i) +a_0 \end{equation} where
$a_i=\pi^* (a_i^0), b_i=\pi^* (b_i^0)$ are lift-ups of regular
functions $a_i^0,b_i^0$ on $X$ via the natural projection $\pi : Y
\to X$ (as we mentioned by abusing terminology we shall say that
$a_i$ and $b_i$ themselves are regular functions on $X$).

Let $\Omega = du \wedge dv \wedge \omega$, i.e. it is a volume
form on $Y$. By Lemma \ref{Divergence} there is a volume form
$\omega'$ on $X'$ such that $\Omega |_{X'} = dP \wedge \omega'$.
Furthermore, for any vector field $\mu$ such that $\mu (P)=0$ and
$\nu' =\mu |_{X'}$ we have $\diver_{\omega'} (\nu' )=
\diver_{\Omega } (\mu )|_{X}$. Note also that any vector field
$\nu$ on $X$ generates a vector field $\kappa$ on $Y$ that
annihilates $u$ and $v$. We shall always denote $\kappa |_{X'}$ by
$ \tilde{\nu}$. It is useful to note for further computations that
$u^i\pi^* (\diver_{\omega} (\nu)) =\diver_{\Omega} (u^i\kappa )$
for every $i \geq 0$. Note also that every algebraic vector field
$\lambda$ on $X'$ can be written uniquely in the form
\begin{equation}\label{Field}
\lambda = \tilde{\mu}_0+ \sum_{i=1}^m (u^i\tilde{\mu}_i^1+ v^i
\tilde{\mu}_i^2) + f_0\partial /\partial u + g_0\partial /
\partial v
\end{equation}
where $\mu_0,\mu_i^j $ are algebraic vector fields on $X$, and
$f_0,g_0$ are regular functions on $X'$.

For any algebraic manifold $Z$ with a volume form $\omega$ we
denote by $\LieA (Z)$ (resp. $\LieAO (Z)$) the Lie algebra
generated by algebraic globally integrable vector fields on $Z$
(resp. that annihilates $\omega$) and by $\VFA (Z)$ we denote the
Lie algebra of all algebraic vector fields on $Z$. We have a
linear map $$\widetilde{\Pr} : \VFA (X') \to \VFA (X)$$ defined by
$\widetilde{\Pr} (\lambda ) = \mu_0$ where $\lambda$ and $\mu_0$
are from formula (\ref{Field}).  As it was mentioned in \cite{KaKu2} the following facts  are
straightforward calculations that follow easily from Lemma
\ref{Divergence}.}
\end{nothing*}

\begin{lemma}\label{Facts}
Let $\nu_1, \nu_2$ be vector fields on $X$,  and $f$ be a regular
function on $X$. For $i\geq 0$ consider the algebraic vector
fields
$$\nu_1'=u^{i+1}\tilde{\nu}_1 +u^i\nu_1(p) \partial /\partial v,
\, \, \, \, \nu_2'=v^{i+1} \tilde{\nu}_2 +v^i \nu_2(p) \partial
/\partial u$$ and $\mu_f =f( u\partial /\partial u - v \partial
/\partial v)$ on $Y$. Then

{\rm (i)} $\nu_i'$ and $\mu_f$ are tangent to $X'$ (actually they
are tangent to fibers of $P=uv-p(x)$), i.e., they can be viewed as
vector fields on $X'$;

{\rm (ii)} $\mu_f$ is always globally integrable on $X'$, and
$\nu_i'$ is globally integrable on $X'$ if $\nu_i$ is globally
integrable on $X$;

{\rm (iii)} $\diver_{\omega'} (\mu_f )=0$, $\diver_{\omega'}
(\nu_1')=u^{i+1} \diver_{\omega} (\nu_1)$, $\diver_{\omega'}
(\nu_2')=v^{i+1} \diver_{\omega} (\nu_2)$, and $$\diver_{\omega'}
([\mu_f, \nu_1'] =(i+1)u^{i+1}f \diver_{\omega} (\nu_1), \,
\diver_{\omega'} ([\nu_2', \mu_f ] =(i+1) v^{i+1}f \diver_{\omega}
(\nu_2);$$

{\rm (iv)} we have the following Lie brackets $$[\mu_f , \nu_1']=
(i+1)u^{i+1}f \tilde{\nu}_1 +\alpha_1 \partial /\partial u +
\beta_1 \partial /\partial v,$$ $$[\nu_2', \mu_f]= (i+1)v^{i+1}f
\tilde{\nu}_2 +\alpha_2 \partial /\partial u + \beta_2
\partial /\partial v,$$ where $\alpha_i$ and $\beta_i$ are some regular
functions on $X'$;

{\rm (v)} more precisely, if $i=0$ in formulas for $\nu_1'$ and
$\nu_2'$ then $$[\mu_f , \nu_1']=fu \tilde{\nu}_1-u^2 \nu_1(f)
\partial /\partial u +  \nu_1(fp) \partial /\partial v,$$
$$[ \nu_2', \mu_f]=fv \tilde{\nu}_2-v^2 \nu_2(f) \partial /\partial
v +  \nu_2(fp) \partial /\partial u;$$  and

\begin{equation}\label{TripleBracket}
\widetilde{\Pr} ([[ \mu_f , \nu_1'], \nu_2']) =
\nu_1(fp)\nu_2-\nu_2(fp) \nu_1 +fp [\nu_1, \nu_2].
\end{equation}
\end{lemma}

\section{The proof of Theorem \ref{one}.}
\label{Theorem1}

\begin{nothing*}\label{Notation} {\bf Additional Notation.}
{\rm For every affine algebraic manifold $Z$ let $\C [Z]$ be the
algebra of its regular functions, $\IVFA (Z)$ be the set of
completely integrable algebraic vector fields on $Z$. If there is a
volume form $\omega$ on $Z$ then we denote by $\Div_Z : \VFA (Z)
\to \C [Z]$ the map that assigns to each vector field its
divergence with respect to $\omega$, and set $\IVFAO (Z) = {\rm
Ker} \, \Div_Z \cap \IVFA (Z)$, $\VFAO (Z) = {\rm Ker} \, \Div_Z
\cap \VFA (Z)$. For a closed submanifold $C$ of $Z$ denote by
$\VFAO (Z,C)$ the Lie algebra of vector fields  of zero divergence
on $Z$ that are tangent to $C$. Formula (\ref{StandardForm})
yields a monomorphism of vector spaces $\iota : \C [X']
\hookrightarrow \C [Y]$ and the natural embedding $X
\hookrightarrow X \times (0,0) \subset Y$ generates a projection
${\rm Pr} : \C [Y] \to \C [X]$. Note that $\Pr (\iota (f))=a_0$ in
the notation of formula (\ref{StandardForm}).}

\end{nothing*}

\begin{lemma}\label{add.2} Let $\lambda$ be a vector field on
$X'\subset X \times \C_{u,v}^2$ given by formula (\ref{Field}).
Suppose that $\omega_0$ is a volume form on $X$ and a volume
 form  $\omega$ on $X'$
coincides with the pull-back of the volume form $\omega_1:=(\omega_0
\wedge du)/u$ on $Z:=X \times \C_u$ under the natural projection
(i.e. $\omega$ constructed as in Remark \ref{add.1}). Then $\Pr
(\diver_{\omega} (\lambda ))=\diver_{\omega_0} (\mu_0 )$. In
particular, if $\diver_{\omega} \lambda =0$ then
$\diver_{\omega_0} (\mu_0 )=0$.

\end{lemma}

\begin{proof} The natural projection $\sigma : X' \to Z$ is an
affine modification whose restriction over $X \times \C_u^*$ is an
isomorphism. Hence $\lambda$ is the pull-back of the following
vector field $$ \kappa  = \tilde{\mu}_0+ \sum_{i=1}^m
u^i(\tilde{\mu}_i^1+ \tilde{\mu}_i^2/p^i) + f_0\partial /\partial
u$$ on $Z$. Thus it suffices to show that $\diver_{\omega_0}
(\mu_0 )= {T_0} (\diver_{\omega_1} (\kappa ))$ where ${T_0} :
 \C (X )[u,u^{-1}] \to \C (X)$ assigns to each Laurent
polynomial in $u$ its constant term. By (\ref{add.3})
$\diver_{\omega_1} (\kappa )=\diver_{\omega_0\wedge du} (\kappa
)-\kappa(u)/u=\diver_{\omega_0\wedge du} (\kappa )-f_0/u $. Hence
$${T_0}(\diver_{\omega_1} (\kappa ))=\diver_{\omega_0} (\mu_0
)+{ T_0}(\partial f_0/\partial u) -{T_0} (f_0/u).$$ The
desired conclusion follows now from the obvious fact that
${T_0}(\partial f_0/\partial u)={T_0} (f_0/u)$.

\end{proof}

\begin{proposition}\label{A} Let $C$ be the smooth zero locus of
$p$ in $X$. Suppose also that the following conditions hold:

{\rm (A1)} the linear space $\VFAO (X,C)$ is generated by vector
fields that are of the form $ \widetilde{\Pr} ([[\mu_f , \nu_1'],
\nu_2'])$ where $\mu_f$ and $\nu_i'$ are as in formula
(\ref{TripleBracket}) from Lemma \ref{Facts} with $\nu_i \in
\IVFAO (X)$;

{\rm (A2)} $\VFA (X)$ is generated by $\IVFAO (X)$ as a module
over $\C [X]$.

Then $\LieAO (X')$ coincides with ${\VFAO} (X')$, i.e., $X'$ has
the algebraic volume density property.

\end{proposition}

\begin{proof}
Let $\lambda , f_0$, and $g_0$ be as in formula (\ref{Field}) and
$\Lambda =\iota (\lambda) $ be the extension of $\lambda$ to $Y$
also given by formula (\ref{Field}). By formula
(\ref{StandardForm}) $f_0$ and $g_0$ can be written uniquely in
the form
$$ f_0=\sum_{i=1}^m (a_iu^i+b_iv^i) +a_0 \, \, {\rm and} \, \,
g_0=\sum_{i=1}^m ( \hat{a}_iu^i+ \hat{b}_iv^i) +\hat{a}_0$$ where
$a_i, \hat{a}_i, b_i, \hat{b}_i \in \C [X]$.

Since $\Lambda$ is a vector field tangent to $X'=P^{-1}(0)$ we
have $\Lambda (P)|_{X'}=0$. Thus $0=\Pr (\iota (\Lambda
(P)|_{X'}))=p (a_1 + \hat{b}_1)-\mu_0 (p) $ (recall that
$P=uv-p(x)$). Hence $\mu_0(p)$ vanishes on $C$, i.e. by Lemma
\ref{add.2} $\mu_0 \in \VFAO (X,C)$. Let $\mu_f, \nu_i'\in \IVFAO
(X')$ be as in Lemma \ref{Facts}. Condition (A1) implies that
adding elements of the form $[[\mu_f ,\nu_1'],\nu_2'] $ to
$\lambda$ we can suppose that $\mu_0=0$.
Using condition (A2) and Lemma \ref{Facts} (iv) we can make
$\mu_i^j=0$ by adding fields of the form $[\mu_f , \nu_i']$ with
$\nu_i \in \IVFAO (X)$. Note that this addition leaves not only
$\mu_0$ equal to 0 but also $\diver_{\omega'} (\lambda )$ equal to
0, since $ \diver_{\omega'} ([\mu_f , \nu_i'])=0$ as soon as
$\diver_{\omega} (\nu_i)=0$. Hence $\lambda = f\partial /\partial
u + g\partial /
\partial v $ and $\Lambda (P)|_{X'} =fv+gu=0$.

Using formula (\ref{StandardForm}) one can see that $f$ must be
divisible by $u$, and $g$ by $v$. That is, there exists a regular
function $h$ on $X'$ for which $f=uh$ and $g=-vh$. Hence $\lambda=
h(u\partial /\partial u -v
\partial /\partial v)$. Note that $\Lambda (P)=0$ now. Thus
by Lemma \ref{Divergence}
$$0=\diver_{\omega'} (\lambda )= \diver_{\Omega} (\Lambda )|_{X'}=
(u {\frac {\partial h}{\partial u}} -v{\frac {\partial h}{\partial v}})|_{X'}.$$ Taking
$h$ as in formula (\ref{StandardForm}) we see that $h$ is
independent of $u$ and $v$. Thus $\lambda$ is integrable and of
zero divergence by Lemma \ref{Facts} (ii)-(iii).

\end{proof}

Taking commutative vector fields $\nu_1$ and $\nu_2$ from $\IVFAO
(X)$ in formula (\ref{TripleBracket}) we have the following.

\begin{lemma}\label{A'}
Condition (A1) in Proposition \ref{A}   holds if $\VFAO (X,C)$ is
generated as a linear space by vector fields of the  form
$\nu_1(fp)\nu_2-\nu_2(fp) \nu_1$ where the vector fields $\nu_1,
\nu_2 \in \IVFAO (X)$ are commutative.
\end{lemma}

It is more convenient for us to reformulate this new condition in
terms of differential forms for which we need some extra facts.
Let $\iota_{\nu}$ be the inner product with a vector field $\nu$ on
$X$. Recall the following relations between the outer
differentiation ${\rm d}$, the Lie derivative $L_{\nu}$ and
$\iota_{\nu}$  \begin{equation}\label{add.4a} L_{\nu}= {\rm
d}\circ \iota_{\nu} + \iota_{\nu} \circ {\rm  d} \, \, \, {\rm
and} \, \, \, [L_{\nu_1}, \iota_{\nu_2}]=\iota_{[\nu_1,
\nu_2]}.\end{equation} Set $\omega_{\nu} =\iota_{\nu} (\omega )$.
Then by formula (\ref{DefDivergence}) we have $$\diver_{\omega}
(\nu )\omega ={\rm d}\circ \iota_{\nu}( \omega ) + \iota_{\nu}
\circ {\rm d} (\omega )={\rm d}( \omega_{\nu} ).$$ Thus we have
the first statement of the following.

\begin{lemma}\label{add.5}\hfill

{\rm (1)} A vector field $\nu$ is of zero divergence if and only
if the form $\omega_{\nu}$ is closed.

{\rm (2)} Furthermore, for a zero divergence field $\nu$
and every regular function $f$ on $X$ we have
${\rm d} ( \omega_{f\nu} ) = \nu (f) \omega$.

{\rm (3)} Let $\nu_1, \nu_2 \in \IVFAO (X)$ be commutative and
$\kappa=\nu_1(fp)\nu_2-\nu_2(fp) \nu_1$. Then ${\rm d} (
\iota_{\nu_1} \circ \iota_{\nu_2} (fp\omega))=\iota_{\kappa}
(\omega )$ where $p, f \in \C [X]$.

\end{lemma}

\begin{proof} Indeed, by (\ref{add.4a})
$${\rm d} ( \omega_{f\nu} ) = L_{f\nu} (\omega )- \iota_{f \nu} {\rm d} (\omega)=
L_{f\nu} (\omega )=\diver_{\omega} (f\nu )\omega = (f\diver_{\omega} (\nu )+\nu (f) ) \omega =\nu (f) \omega$$
which is (2).

Again by (\ref{add.4a}) we have $${\rm d} \circ
\iota_{\nu_1} \circ \iota_{\nu_2} (fp\omega)=L_{\nu_1} \circ
\iota_{\nu_2} (fp\omega) - \iota_{\nu_1} \circ {\rm d} \circ
\iota_{\nu_2} (fp\omega).$$ Then $$L_{\nu_1} \circ \iota_{\nu_2}
(fp\omega)=L_{\nu_1} (fp) \iota_{\nu_2} (\omega) +fp L_{\nu_1}
\circ \iota_{\nu_2} (\omega)$$ and $$ L_{\nu_1} \circ
\iota_{\nu_2} (\omega)= \iota_{\nu_2}
L_{\nu_1}(\omega)+\iota_{[\nu_1,\nu_2]}(\omega )=0$$ since
$[\nu_1,\nu_2]=0$ and $L_{\nu_i}(\omega)=0$. Similarly
$$\iota_{\nu_1} \circ {\rm d} \circ
\iota_{\nu_2} (fp\omega))=L_{\nu_2}(fp) \iota_{\nu_1} (\omega)+
fp\iota_{\nu_1} \circ L_{\nu_2} (\omega )-\iota_{\nu_1}  \circ
\iota_{\nu_2} \circ {\rm d}(fp\omega)) =L_{\nu_2}(fp)
\iota_{\nu_1} (\omega).$$ Therefore, $${\rm d} \circ \iota_{\nu_1}
\circ \iota_{\nu_2} (fp\omega)=L_{\nu_1} (fp) \iota_{\nu_2}
(\omega)-L_{\nu_2}(fp) \iota_{\nu_1} (\omega)={\nu_1} (fp)
\iota_{\nu_2} (\omega)-{\nu_2}(fp) \iota_{\nu_1} (\omega)$$ which
yields the desired conclusion.

\end{proof}

Suppose that $\Omega^{q}(X)$ is the sheaf of algebraic $q$-forms
on $X$, $\Omega^{q}_i(X)$ is its subsheaf that consists of forms
that vanish on $C$ with multiplicity at least $i$ for $i \geq 1$,
and vanish on all elements $\Lambda^{n-1}TC \subset
\Lambda^{n-1}TX$ for $i=0$ where $\Lambda^{q}TX$ is the $q$-h
wedge-power of $TX$, i.e. the set of $q$-dimensional subspaces of
the tangent bundle. For every sheaf $\cF$ on $X$ denote by
$\Gamma^0(X, \cF )$ the space of global sections. That is,
$\Gamma^0(X,\Omega_1^{n-2}(X))$ is the subset of
$\Gamma^0(X,\Omega^{n-2} (X))$, that consists of forms divisible
by $p$, and $\Gamma^0(X, \Omega_0^{n-1}(X))$ is the set of algebraic
$(n-1)$-forms on $X$ whose restriction to the zero fiber $C$ of
$p$ yields a trivial form on $C$.

As a consequence Lemma \ref{add.5} we have the following fact.

\begin{lemma}\label{add.6} Let  $\kappa_i^f=\nu_1^i(fp)\nu_2^i-\nu_2^i(fp) \nu_1^i$
and let the following condition hold:

\noindent  {\rm (B )} there exists a collection $\{ \nu_1^i,
\nu_2^i \}_{i=1}^m$ of pairs of commutative vector fields from
$\IVFAO (X)$ such that the set $\{ \iota_{\nu_1^i} \circ
\iota_{\nu_2^i} (\omega) \}_{i=1}^m$ generates the space of
algebraic $(n-2)$-forms $\Gamma^0(X,\Omega^{n-2} (X))$ on $X$ as $\C[X]$-module.

Then the image of $\Gamma^0(X,\Omega_1^{n-2}(X))$ under the outer
differentiation ${\rm d} : \Gamma^0(X, \Omega^{n-2}(X)) \to
\Gamma^0(X,\Omega^{n-1}(X) )$ is generated as a vector space by
($n-1$)-forms $\{ \iota_{\kappa_i^f} (\omega ) \}_{i=1}^n, f \in \C[X]$.

\end{lemma}

\subsection{Application of Grothendieck's theorem}
Let $\cZ^0(X, \Omega_0^{n-1}(X))$ be the subspace of closed algebraic
$(n-1)$-forms in $\Gamma^0(X, \Omega^{n-1}_0)$.
Clearly, for every algebraic vector
field $\nu \in \IVF_{\omega} (X)$ tangent to $C$ we have
$\omega_{\nu} \in \cZ^0(X,\Omega_0^{n-1}(X))$. Our aim now is to
show that under mild assumption the homomorphism $${\rm d} :
\Gamma^0(X,\Omega_1^{n-2}(X)) \to \cZ^0(X,\Omega_0^{n-1}(X)) $$ is
surjective and, therefore, condition (A1) from Proposition \ref{A}
follows from condition (B) from Lemma \ref{add.6}.
Denote by $\cF'_i$ (resp. $\cF_i$) the space
of algebraic sections of $\Omega_{n-1-i}^{i}$ (resp. $\Omega^i$)
over $X$. Note that the outer differentiation ${\rm d}$ makes
$$\cF' (*) := \ldots  \to \cF_i' \to \cF_{i+1}'\to \ldots \, \, \, {\rm and}
\, \, \, \cF (*) := \ldots \to \cF_i \to \cF_{i+1} \to \ldots $$
complexes, and that the surjectivity we need for condition (B) would follow
from $H^{n-1} (\cF' (*))= 0$.

\begin{proposition}\label{sur.1} Let $H^{n-1}(X, \C ) =0$ and let  the
homomorphism $H^{n-2}(X,\C ) \to  H^{n-2}(C,\C )$ generated by the
natural embedding $C \hookrightarrow X$ be surjective. Then
$H^{n-1} (\cF' (*))$   $=0$.
\end{proposition}

\begin{proof}
Consider the following short exact sequence of complexes $0 \to
\cF'(*) \to \cF (*) \to \cF''(*) \to 0$ where
$\cF_i''=\cF_i/\cF_i'$ in complex $\cF'' (*)$. This implies the
following long exact sequence in cohomology
$$\ldots \to H^{n-2}(\cF (*) ) \to H^{n-2}(\cF'' (*) ) \to H^{n-1}(\cF' (*)
)\to H^{n-1}(\cF (*) ) \to \ldots , $$ i.e. we need to show (i)
that the homomorphism $H^{n-2}(\cF (*) ) \to H^{n-2}(\cF'' (*) )$
is surjective and (ii) that $ H^{n-1}(\cF (*) )=0$. By the
Grothendieck theorem \cite{Gro} De Rham cohomology on smooth affine
algebraic varieties can be computed via the complex of algebraic
differential forms, i.e. $ H^{n-1}(\cF (*) )=H^{n-1}(X,\C )$ which
implies (ii). Similarly, $ H^{n-2}(\cF (*) )=H^{n-2}(X,\C )$. Note
that $\cF''_i= \cF_i/ (p^{n-1-i} \cF_i)$ for $i \leq n-2$. In
particular, modulo the space $\cS$ of (the restrictions to $C$ of)
algebraic $(n-2)$-form that vanish on $\Lambda^{n-2}TC$ the term
$\cF''_{n-2}$ coincides with the space $\cT$ of algebraic
$(n-2)$-forms on $C$ (more precisely, we have the following exact
sequence $0\to \cS \to \cF''_{n-2} \to \cT\to 0$). One can see that
each closed $\tau \in \cS$ is of form $ {\rm d}p \wedge \tau_0$
where $\tau_0$ is a closed $(n-3)$-form on $C$. Hence $\tau ={\rm
d} (p \tau_0 )$ (where $p\tau_0$ can be viewed as an element of
$\cF''_{n-3}=\cF_{n-3}/ (p^2 \cF_{n-3})$) is an exact form. Thus
the $(n-2)$-cohomology of complex $\cF''(*)$ coincides with the
$(n-2)$-cohomology of the algebraic De Rham complex on $C$ and,
therefore, is equal to $H^{n-2} (C, \C )$ by the Grothendieck
theorem. Now homomorphism from (i) becomes $H^{n-2}(X,\C ) \to
H^{n-2}(C,\C )$ which implies the desired conclusion.

\end{proof}

Thus we have Theorem \ref{one} from Introduction as a consequence  of
the following more general fact (which gives, in particular,
an affirmative answer  to an open question of \textsc{Varolin} (\cite{V3}, section 7) who asked whether the hypersurface
$\{(a, b, c, d) \in \C^4 : a^2c-bd = 1 \}$ in $\C^4$ has the volume density property).

\begin{theorem}\label{hypersurface}
Let $X$ be an $n$-dimensional smooth affine algebraic variety with
$H^n(X,\C)=0$ and a volume form $\omega$ satisfying
conditions

{\rm (B )} there exists a collection $\{ \nu_1^i,
\nu_2^i \}_{i=1}^m$ of pairs of commutative vector fields from
$\IVFAO (X)$ such that the set $\{ \iota_{\nu_1^i} \circ
\iota_{\nu_2^i} (\omega) \}_{i=1}^m$ generates the space of
algebraic $(n-2)$-forms $\Gamma^0(X,\Omega^{n-2} (X))$ on $X$ as $\C[X]$-module;

{\rm (A2)} $\VFA (X)$ is generated by $\IVFAO (X)$ as a module
over $\C [X]$.\footnote{
Clearly, the standard volume form on $\C^n$ satisfies both these conditions.}

Suppose also that $p$
is a regular function on $X$ with a smooth reduced zero fiber $C$ such that the
homomorphism $H^{n-2}(X,\C ) \to  H^{n-2}(C,\C )$ generated by the
natural embedding $C \hookrightarrow X$ is surjective.  Let $X' \subset X \times \C^2_{u,v}$
be the hypersurface given by $uv=p$ and let $\omega'$ be the pullback of the form
$\omega \wedge du /u$ on $Z =X \times \C_u$ under the natural projection $X' \to Z$ \footnote{One can
check that $\omega'\wedge dP =\omega \wedge du \wedge dv|_{X'}$ where $P=uv -p$.}. Then
$X'$ has the algebraic $\omega'$-density property.
\end{theorem}

\subsection{Algebraic volume density for $SL_2(\C )$ and $PSL_2(\C )$.}
Since $X'=SL_2(\C )$ is isomorphic to the hypersurface $uv=x_1x_2+1 =: p (\bar x )$ in $Y=\C_{\bar x , u,v}^4=X\times \C_{u,v}^2$ with $\bar x=(x_1,x_2)$ and $X=\C^2_{\bar x}$, Theorem \ref{one} implies that $SL_2(\C )$ has the algebraic volume density property with respect to the volume form $\omega'$ on $X'$ such that $\omega' \wedge {\rm d} P= \Omega$ where $P=uv- p (\bar x)$ and $\Omega = {\rm d} x_1 \wedge {\rm d} x_2 \wedge {\rm d} u \wedge {\rm d} v$ is the standard volume form on $\C^4$. On the other hand by Remark \ref{add.1} (1) we can consider forms $({\rm d} x_1 \wedge {\rm d} x_2 \wedge {\rm d} u)/u$,  $({\rm d} x_1 \wedge {\rm d} x_2 \wedge {\rm d} v)/v$, etc.. Each of these forms coincides with $\omega'$  up to a sign because their wedge-products with ${\rm d} P$ are $\pm \Omega$. Note that $({\rm d} x_1 \wedge {\rm d} x_2 \wedge {\rm d} u)/u$ is invariant with respect to the $\C_+$-action on $SL_2(\C )$ given by $(x_1,x_2,u,v) \to (x_1, x_2+tx_1, u, v+tu), \, t \in \C_+$ which is generated by multiplications of a $\C_+$-subgroup of $SL_2(\C )$.  Thus $\omega'$ is invariant with respect to such multiplications. Similarly, consideration of $({\rm d} x_1 \wedge {\rm d} x_2 \wedge {\rm d} v)/v$ yields invariance with respect to the $\C_+$-action $(x_1,x_2,u,v) \to (x_1+tx_2,x_2,u+tv,v)$, etc.. This implies that $\omega'$ is invariant with respect to multiplication by any element of $SL_2(\C )$ and we proved the following result, which is originally due to Varolin (\cite{V3}, Theorem 2).

\begin{proposition}\label{sl2'} Group $SL_2(\C )$ has the algebraic volume density property with respect to the  invariant volume form.

\end{proposition}

Furthermore, since the commutative  vector fields $\nu_1= \partial /\partial x_1$ and $\nu_2 = \partial / \partial x_2$ on $X=\C_{\bar x}^2$ satisfy condition (B) of Lemma \ref{add.6} we see that any vector field $\mu_0$ tangent to the zero fiber $C$ of $p$ is of form $\nu_1(fp)\nu_2 -\nu_2 (fp)\nu_1$ where $f$ is a polynomial on $X$. This fact will used in the next unpleasant computation which is similar to  the argument in Proposition \ref{A}.

\begin{proposition}\label{psl2} Group $PSL_2(\C )$ has the algebraic volume density property with respect to the  invariant  volume form.

\end{proposition}

\begin{proof}
 Consider now $X'' =X'/\Z_2 \simeq PSL_2(\C )$ where the $\Z_2$-action on
$X'$ given by $(u,v,\bar x )\to (-u, -v,-\bar x )$.
Note that  $\C [X'']$ can be viewed as the subring of $\C [X']$ generated by monomials of even degrees. Hence completely integrable vector  fields of form
$$ \nu_1'=u^{i+1}\partial /\partial x_k +u^i {\frac{\partial p}{\partial x_k}} \partial / \partial v
\, \, \, {\rm  and} \, \, \,
\nu_2'=v^{i+1}\partial /\partial x_j + v^i{\frac{\partial p}{\partial x_j}} \partial / \partial u$$
$${\rm (resp.} \, \, \,  \nu_1''=u^{i+1}x_j\partial /\partial x_k +u^i x_j{\frac{\partial p}{\partial x_k}} \partial / \partial v\, \, \, {\rm  and}\, \, \,
 \nu_2''=v^{i+1}x_k\partial /\partial x_j +v^i x_k{\frac{\partial p}{\partial x_j}} \partial / \partial u
 {\rm )} $$ on $X'$ with even (resp. odd) $i$ can viewed as fields on $X''$. The same is true for $\mu_f$ from Lemma \ref{Facts} provided $f$ is a linear combination of monomials of even degrees.  Fields $\nu_1', \nu_2', \mu_f$ are of zero divergence. If $j \ne k$ the same holds for  $\nu_1''$ and $\nu_2''$ .  Any algebraic vector field $\lambda$ on $X''$ can be viewed as a vector field on $X'$ and, therefore, it is given by formula (\ref{Field}).  Since this field on $X'$ came from $X''$ each $\tilde \mu_i^k$ (resp. $\tilde \mu_0$)  in that formula consists of summands of form $q(\bar x) \partial / \partial x_k$ where polynomial $u^iq(\bar x)$  (resp. $v^iq(\bar x)$) is a linear combination of monomials of odd degrees.  Our plan is to simplify the form of a zero divergence field $\lambda$ on $X''$ by adding elements of the Lie algebra generated by  fields like $\mu_f, \nu_1', \nu_2', \nu_1'', \nu_2''$.

Recall that $\tilde \mu_0$ is generated by a field $\mu_0$ on $X$ and it was shown in the proof of Proposition \ref{A} that  $\mu_0$ is tangent to $C$. Hence, as we mentioned before $ \mu_0=\nu_1(fp)\nu_2 -\nu_2 (fp)\nu_1$. Furthermore, if $\lambda$ comes from a field on $X''$ polynomial $f \in \C [x_1,x_2]$ must contain monomials of even degrees only.  Thus by virtue of Lemma \ref{Facts} (v) adding to $\lambda$ vector fields of form $[[\mu_f , \nu_1'], \nu_2']$ we can suppose that $\tilde \mu_0=0$ without changing the divergence of $\lambda$.

Following the pattern of the proof of Proposition \ref{A} let us add to $\lambda$ the zero divergence fields of form $[\mu_f, \nu_l']$ and $[\mu_f, \nu_l'']$. Since we have to require that $j\ne k$ in the definition of $\nu_1''$ and $\nu_2''$ we cannot eliminate summands $u^{i+1}\tilde \mu_{i+1}^1$ and $v^{i+1} \tilde \mu_{i+1}^2$ completely. However, Lemma \ref{Facts} (iv)  shows that  after such addition one can suppose that  $u^{i+1}\tilde \mu_{i+1}^1$ vanishes for even $i$, and for odd $i$ it is a linear combination of terms of form $ u^{i+1}x_k^m\partial /\partial x_k$ where $m$ is odd (and similarly for  $v^{i+1}\tilde \mu_{i+1}^2$).

Consider the semi-simple vector field $\nu =x_1\partial /\partial x_1 -x_2\partial /\partial x_2$ on $X$. Then $\nu' =u^{i+1}\tilde  \nu + u^i\nu (p) \partial /\partial v$ is a completely integrable zero divergence vector field on $X'$ and for odd  $i$ it can be viewed as a field on $X''$. Set $f=x_1^{m-1}$. By Lemma \ref{Facts} (iv) $$[\mu_f, \nu']= (i+1)u^{i+1}x_1^{m-1} \tilde \nu + \alpha \partial / \partial u + \beta \partial /\partial v. $$Thus adding a multiple of $[\mu_f ,\nu' ]$ to $\lambda$ we can replace terms $u^{i+1}x_1^m\partial /\partial x_1$ in $u^{i+1} \tilde \mu_{i+1}^1$ by  $u^{i+1}x_1^{m-1}x_2\partial /\partial x_2$. If $m\geq 2$ the latter can be taken care of by adding fields of form $[\mu_f , \nu_i'']$. If $m=1$ we cannot eliminate immediately terms like $u^ix_1 \partial /\partial x_1$ or $u^ix_2 \partial / \partial x_2$, but adding fields of form $cu^i \tilde\nu$ where $c$ is a constant we can suppose that only one of these terms is present. The same is true for  similar terms with $u$ replaced by $v$. Thus adding elements from $\Lie^{\omega} (X'')$  we  can reduce $\lambda$ to a zero divergence  field of the following form
$$\lambda =\sum_{i\geq 1} (c_i u^ix_1\partial /\partial x_1 +d_iv^ix_2\partial /\partial x_2) +g_1 \partial /\partial u + g_2 \partial /\partial v$$ where constants $c_i$ and $d_i$ may be different from zero only for even indices $i$ and  by formula (\ref{StandardForm}) $g_k=\sum_{i\geq 1}(a_i^k(\bar x) u^i+b_i^k(\bar x) v^i)   +a_0^k(\bar x)$ with $a_i^k$ and $b^k_i$ being polynomials on $X$.
Since divergence $\diver_{\omega} \lambda =0$ we immediately have $a_{i+1}^1=-c_i/(i+1)$ and $b_{i+1}^2=-d_i/(i+1)$, i.e. these polynomials are constants.

Consider now an automorphism of $X''$ (and, therefore, of $X'$) given by $(u,v,x_1,x_2) \to (-x_1,x_2,-u,v)$, i.e. it exchanges the role of pairs $(u,v)$ and $(x_1,x_2)$. It transforms $\lambda$ into a field $$\sum_{i\geq 1} (a_{i+1}^1 x_1^{i+1}\partial /\partial x_1 +b_{i+1}^2x_2^{i+1} \partial /\partial x_2) + \lambda_0$$  where $\lambda_0$ does  not contain nonzero summands of form $a x_1\partial /\partial x_1$ (resp. $b x_2\partial /\partial x_2$) with $a$ (resp. $b$)  being a regular function on $X''$ non-divisible by $x_2$ (resp. $x_1$). Hence adding fields of form $[[\mu_f, \nu_1'],\nu_2']$, $[\mu_f, \nu_k']$, and $[\mu_f, \nu_k'']$ as before we can suppose that $\tilde \mu_0$ and each $\tilde \mu_i^1$ and $\tilde \mu_i^2$ are equal to zero, i.e. $\lambda = e \partial /\partial u + g\partial /\partial v$. Furthermore, arguing as in Proposition \ref{A} we see $\lambda = h (u \partial /\partial u - v \partial /\partial v)$ where $h$ is a polynomial on $X$, i.e. $\lambda$ is completely integrable. Since by construction it is a vector field on $X''$, we have proved that $X''$ possesses  the algebraic volume density property.

\end{proof}

\section{Two basic  facts about the algebraic volume density property}\label{twofacts}

By considering differential forms and vector fields in local
coordinate systems one can see that the map $\nu \to \omega_{\nu}
:=\iota_{\nu} (\omega )$ is bijective and, therefore, establishes
a duality between algebraic (resp. holomorphic) vector fields and
the similar $(n-1)$-forms on $X$. This duality in combination with
the Grothendieck theorem \cite{Gro} enables us to prove another important fact.

\begin{proposition}\label{approx}
For an affine algebraic manifold $X$ equipped with an algebraic
volume form $\omega$ the algebraic volume density property implies the
volume density property (in the holomorphic sense).
\end{proposition}

\begin{proof} We need to show that any holomorphic vector
field $\mu$ such that $\mu (\omega )=0$ can be approximated by an
algebraic vector field $\nu$ with $\nu ( \omega )=0$. Since the
form $\omega_{\mu}$ is closed, by the Grothendieck theorem one can
find a closed algebraic $(n-1)$-form $\tau_{n-1}$ such that
$\omega_{\mu}-\tau_{n-1}$ is exact, i.e. $\omega_{\mu}-\tau_{n-1}=
{\rm d} \, \tau_{n-2}$ for some holomorphic $(n-2)$-form
$\tau_{n-2}$. Then we can approximate $\tau_{n-2}$ by an algebraic
$(n-2)$-form $\tau_{n-2}'$. Hence the closed algebraic
$(n-1)$-form $\tau_{n-1}+ {\rm d} \, \tau_{n-2}'$ yields an
approximation of $\mu$. By duality $\tau_{n-1}+ {\rm d} \,
\tau_{n-2}'$ is of form $\omega_{\nu}$ for some algebraic vector
field $\nu$ (approximating $\mu$ and by Lemma \ref{add.5} (1) $\nu$ is of zero
$\omega$-divergence which is the desired conclusion.

\end{proof}

\begin{lemma} \label{span}
If $X$ has the algebraic volume density property, then there exist finitely many
algebraic vector fields $\sigma_1, \ldots , \sigma_m \in \Lie_{alg}^{\omega} (X)$
that generate $\VFA (X)$  as a $\C [X]$-module.
\end{lemma}

\begin{proof}
Let $n = \dim X$. We start with the following.

{\em Claim.} The space of  algebraic  fields of zero divergence
generates the tangent space of $X$ at each point.

Let $x \in X$ and $U$ be a Runge neighborhood of $x$ such that $H^{n-1} (U, \C) = 0$ (for example
take a small sublevel set of a strictly plurisubharmonic  exhaustion
function on $X$ with minimum at $x$).
Shrinking $U$ we can assume that in some holomorphic coordinate system $z_1, \ldots ,z_n$
on $U$
the form $\omega \vert_U$ is the standard volume  ${\rm d} z_1\wedge  \ldots \wedge {\rm d} z_n$.
Thus the holomorphic vector fields $\partial /\partial z_i$ on $U$ are of zero divergence
and they span the tangent space at $x$. We need to approximate them by global algebraic
fields of zero divergence on $X$ which would yield our claim. For that let $\nu \in \VF_{hol}^\omega (U)$. The inner product $\iota_\nu (\omega ) =: \alpha$ is by Lemma \ref{add.5} (1) a closed
$(n-1)$-form  on $U$ and since $H^{n-1} (U, \C) = 0$ we  can
find an $(n-2)$-form $\beta$ on $U$ with ${\rm d} \beta = \alpha$.
Since $U$ is Runge in $X$ we can also
approximate $\beta$ by a global algebraic $(n-2)$-form $\tilde \beta$ (uniformly on compacts in $U$). Then the closed algebraic $(n-1)$-form ${\rm d} \tilde\beta$
approximates $\alpha$ and the unique algebraic vector field $\theta$ defined by
$\iota_\theta (\omega ) = {\rm d} \tilde \beta$ approximates $\nu$. Since
${\rm d} \tilde \beta$ is closed, the field $\theta$ is of zero divergence which
concludes the proof of the Claim.

Now it follows from the Claim and  the algebraic volume density property
that  there are $n$ vector fields in $ \Lie_{alg}^{\omega} (X)$ which span the tangent space at a given
point $x \in X$. By standard induction on the dimension,
adding more fields to span the tangent spaces at points where it was not spanned yet,
we get the assertion of the lemma.

\end{proof}

\begin{remark}\label{08.09a} The similar fact holds for the (holomorphic) density
property, because the same argument implies the holomorphic version of the Claim. In fact,
we can say more: completely integrable holomorphic fields of zero divergence on a Stein manifold
$X$ with the density property generate the tangent space at each point.
Indeed, since  by the Claim
such a tangent space is generated by elements of $\LieHO (X)$ it suffices to show
that every Lie bracket $[\nu , \mu ]$ of completely holomorphic integrable vector fields
of zero divergence
can be approximated by a linear combination of such fields which follows immediately from
the equality $[ \nu , \mu ] = \lim_{t \to 0}  {\frac{\phi_t^* (\nu ) - \nu}{t}}$ where $\phi_t$ is
the phase flow generated by $\mu$. The reason why we cannot make this stronger version of
the Claim in the algebraic category is that for an algebraic vector field $\mu$ the phase
flow $\phi_t$ may be only holomorphic. Another interesting fact is that a set of completely integrable holomorphic
zero divergence vector fields which is needed for generation of tangent spaces at each point of $X$ can be chosen finite.

\end{remark}

Let us suppose that $X$ and $Y$ are affine algebraic
manifolds equipped with
volume forms  $\omega_X$ and $\omega_Y$ respectively.

\begin{proposition}\label{pro1} Suppose that  $X$ (resp. $Y$) has the algebraic $\omega_X$ (resp. $\omega_Y$) volume density property. Let $\omega =\omega_X \times \omega_Y$. Then $X \times Y$ has the
algebraic volume density property relative to $\omega$.

\end{proposition}

\begin{proof}  By Lemma \ref{span} we can suppose that $\sigma_1, \ldots , \sigma_m \in \Lie_{alg}^{\omega_X} (X)$
(resp.  $\delta_1, \ldots , \delta_n \in \Lie_{alg}^{\omega_Y} (Y)$) generate $\VFA (X)$ as a $\C [X]$-module
(resp. $ \VFA (Y)$ as a $\C [Y]$-module).

Denote by $F_Y$ the vector subspace (over $\C$)
 of $\C [Y]$ generated by $ {\rm Im } \, \delta_1, \ldots , {\rm Im } \,  \delta_n$.
Then $\C [Y] = F_Y \oplus V$ where $V$ is another subspace whose
basis is $v_1, v_2, \ldots$. Set $F_Y'=\C [X] \otimes F_Y$ and
$V'=\C [X] \otimes V$, i.e. the algebra of regular functions on $X
\times Y$ is $A=\C [X]\otimes \C [Y] = F_Y' \oplus V'$.  Let $f_i
\in \C[X]$ and $g_j \in \C[Y]$. Note that $f_i$ is in the kernel of all globally integrable fields used in the Lie combination for $\delta_i$ and thus $f_i \delta_i \in \Lie_{alg}^{\omega_Y} (Y)$, analogously
$g_i \sigma_i \in \Lie_{alg}^{\omega_X} (X)$. The fields
$\delta_i$ and $\sigma_j$ generate (vertical and horizontal)
vector fields on $X \times Y$ that are denoted by the same
symbols. Consider $$ [f_i\delta_i, g_j\sigma_j]= \delta_i
(f_ig_j)\sigma_j-\sigma_j (f_ig_j)\delta_i .$$ By construction
$\delta_i$ and $\sigma_j$ are commutative and moreover $\textrm{Span} \, f_i \cdot g_i =\C [X\times Y]$. Hence   the coefficient before
$\sigma_j$ runs over ${\rm Im } \, \delta_i$ and, therefore, for
any $\alpha_1', \ldots , \alpha_n' \in F_Y'$ there are $\beta_1',
\ldots , \beta_m' \in A$ such that the vector field

$$\sum_j \alpha_j' \sigma_j -\sum_i \beta_i' \delta_i$$ belongs to $ \LieAO (X \times Y )$.
Thus adding vector fields of this form to a given vector field
$$\nu =\sum_j \alpha_j \sigma_j -\sum_i \beta_i \delta_i$$  from $\VFAO (X \times Y) $
we can suppose that each $\alpha_j \in V'$.
Hence one can rewrite $\nu$ in the following form $$\nu =\sum_l\sum_j
(h_{jl} \otimes v_l) \sigma_j -\sum_i \beta_i \delta_i$$ where each $h_{jl} \in \C [X]$. Then one has
$$0= \diver \nu = \sum_l (\sum_j
 \sigma_j (h_{jl})) \otimes v_l  -\sum_i \delta_i (\beta_i).$$ Since the first summand is in $V'$ and the last
is in $F_Y'$ we see that $\sum_j
 \sigma_j (h_{jl}) =0$, i.e. each vector field $\sum_j
h_{jl} \sigma_j$ belong to $\VF_{alg}^{\omega_X} (X)$ and by the assumption to $\Lie_{alg}^{\omega_X} (X)$. Hence it suffices to
prove the following

{ \em Claim. } Consider the subspace $ B\subset \VFAO (X \times
Y)$ that consists of vector fields of form $$\nu =\sum_i \beta_i \delta_i\, .$$ Then $B$  is contained in
$\Lie_{alg}^{\omega} (X \times Y)$.

Indeed, consider a closed embedding
of $Y$ into a Euclidean space. Then it generates filtration on $\C
[Y] $ by minimal degrees of extensions of regular functions to
polynomials. In turn this generates filtrations  $B
= \bigcup B_i$ and  $\Lie_{alg}^{\omega_Y} (X\times Y)= \bigcup L_i $. Note
that each $B_i$ or $L_i$ is a finitely generated $\C [X] $-module,
i.e. they generate coherent sheaves on $X$. Furthermore, since $Y$
has algebraic $\omega_Y$-density property we see that the
quotients of $B_i$ and $L_i$ with respect to the maximal ideal
corresponding to any point $x \in X$ coincide. Thus $B_i=L_i$
which implies the desired conclusion.

\end{proof}

Note that up  to a constant factor the completely integrable
vector field $z\partial/\partial z$ on the group $X=\C^*$ is the
only field of zero divergence with respect to the invariant volume
form $\omega = {\frac{{\rm d} z}{z}}$, i.e.,  $X$ has the algebraic volume density
property. Hence we have the following (see also Corollary 4.5 in \cite{V1}).

\begin{proposition}\label{prod.11} For every $n \geq 1$ the torus $(\C^*)^n$ has the
algebraic volume density property with respect to the invariant
form.

\end{proposition}

\section{Algebraic volume density for locally trivial fibrations}

\subsection{Fine manifolds}

If the subspace $V\simeq \C [Y]/F_Y$ in the proof of Proposition \ref{pro1} were trivial so would be the proof  but in  the general case $V \ne 0$.
We shall
see later that Proposition \ref{pro1} can be extended to some locally trivial fibrations with fiber $Y$
for which, in particular,  $ \C [Y]/F_Y$ is at most one-dimensional. Thus we need the following.

\begin{definition}\label{ffi.de1}
Let $Y$ be an affine algebraic manifold with a volume form $\omega$ and $F_Y$ be the subspace of $\C[Y]$ that consists of images of vector fields from $\Lie_{alg}^{\omega}(Y)$.
We say that
$Y$ is a fine manifold if either $\C [Y]=F_Y$ or $\C [Y] $ is naturally
isomorphic to $F_Y\oplus \C$ where
the second summand denotes constant functions on $Y$.

\end{definition}

\begin{lemma}\label{ffi.le1} Let $Y$ be the smooth hypersurface in
$\C_{u,v \bar x}^{n+2}$  given by $P=uv +q(\bar x)-1=0$ where $q(\bar x) =\sum_{i=1}^n
x_i^2$ (i.e. after a coordinate change
$uv +q(\bar x)$ can be replaced by any non-degenerate quadratic form).  Suppose that $Y$  is equipped with a volume $\omega_Y$ such that ${\rm d}P
\wedge \omega =\Omega$ where $\Omega$ is the standard volume form
on $\C^{n}$. Then

{\rm (1)} $Y$ is a fine manifold;

{\rm (2)} $Y/\Z_2$ is a fine manifold where the $\Z_2$-action is given  by $(u,v, \bar x) \to (-u,-v, - \bar x )$.

\end{lemma}

\begin{proof}
Consider the semi-simple vector field $\mu = u\partial /\partial u - v \partial /\partial v$ on
$Y$.
It generates $\Z$-grading of $\C [Y] = \oplus_{i \in \Z} A_i$ such that $\Ker
\mu =A_0$ and $\Image \mu = \oplus_{i \in \Z , i \ne 0} A_i \subset  F_Y$.
Note that $A_0 \simeq \C [x_1, \ldots , x_{n}]$ since $uv =1-q(x) \in
A_0$.  Assume for simplicity that $n\geq 2$ and replace $x_1$ and $x_2$ by $u'=x_1 +\sqrt{-1}x_2$ and $v'=x_1 - \sqrt{-1}x_2$ in our
coordinate system.  Consider the semi-simple vector field  $\mu' = u'\partial /\partial u' - v' \partial /\partial v'$ whose kernel is $A_0'=\C [u,v, x_3, \ldots , x_n]$. Thus monomials containing $u'$ and $v'$ (or, equivalently, $x_1$ or  $x_2$ in the original coordinate system) are in $\Image \mu ' \subset F_Y$.  Repeating this procedure with other $x_i$ and $x_j$ instead of $x_1$ and $x_2$ we see that $F_Y$ contains every nonconstant monomial which is (1).

For (2) note that $\C [Y/\Z_2]$ is the subring of $\C [Y]$ generated by monomials of even degrees and that the semi-simple vector fields that we used preserve the standard degree function. That is, if a monomial $M_1$ of even degree belongs, say, to $\Image \mu$ then $M_1= \mu (M_2)$ where $M_2$ is also a monomial of even degree. This yields (2).

\end{proof}

Since $SL_2(\C )$ is isomorphic to the hypersurface $uv-x_1x_2=1$ in $\C^4_{u,v,x_1,x_2}$ we have the following.

\begin{corollary}\label{ffi.cor1a} Both $SL_2(\C )$ and $PSL_2(\C )$ are fine manifolds.

\end{corollary}

\begin{remark}\label{failed.invariant} In fact for $Y$ equal to $SL_2(\C )$ or $PSL_2(\C )$ we have $\C [Y]/F_Y\simeq \C$.
More precisely, set $F = {\rm Span} \, \{ \nu (f) : \nu \in \VFAO (Y), \, f \in \C [Y] \}$. Note that
vector fields of form $f \nu$
span all algebraic vector fields because of Claim in Lemma \ref{span} (and Lemma \ref{vfi.no0} below).
Therefore, $(n-1)$-forms $\omega_{f\nu}$ generate all algebraic $(n-1)$-forms on $Y$ where $n = \dim Y$.
By Lemma \ref{add.5} (2), ${\rm d} (\omega_{f\nu})= \nu (f) \omega$ which implies that
the image of $\Omega^{n-1}(Y)$ in $\Omega^n (Y)$ under outer differentiation coincides with $F \omega$.
Since ${\rm d} (\Omega^n (Y))=0 $ we have $\C[Y] /F \simeq H^n(Y, \C)$ by the Grothendieck theorem.
By Proposition 4.1 in \cite{KZ} for a smooth hypersurface $Y\subset \C^{m+2}$ given by $uv= p (x) $ we have $H_* (Y) =H_{*-2}(C)$ where $C$ is
the zero fiber of $p$. Thus the universal coefficient formula
implies that  $\dim \C [Y] /F = {\rm rank } \, H^{m-1} (C, \C )$.
For $SL_2 (\C )$ presented as such a hypersurface we have $p(x_1,x_2)=x_1x_2-1$, i.e. $C$ is a hyperbola
and $H^1(C, \C )= \C$ which yields the desired conclusion because $F=F_Y$ for manifolds
with the algebraic volume density property.
\end{remark}

\begin{notation}\label{frank.10}  Further in this section
$X,Y,$ and $W$ are  smooth affine algebraic varieties and  $p : W \to X$
is a locally trivial fibration with fiber $Y$ in the
\'etale  topology. We suppose also that $Y$ is equipped with
a unique up to constant algebraic volume form $\omega_Y$,  and $\VFA (W,p)$ (resp.
$\VF_{alg}^{\omega_Y}(W,p)$) is the space of algebraic vector
fields tangent to the fibers of $p$ (resp. and such that the
restriction to each fiber has zero divergence relative to
$\omega_Y$. ) Similarly $\Lie^{\omega_Y} (W,p)$ will be the Lie algebra
generated by completely integrable vector fields from $\VF_{alg}^{\omega_Y} (W,p)$.
We denote the subspace of $\C [W]$ generated by
functions of form $\{ \Image \nu | \nu \in \Lie^{\omega_Y} (W,p) \}$
by $F(W,p)$.
\end{notation}

\begin{definition}\label{pro2}
A family $\delta_1, \ldots , \delta_n, \ldots  \in \Lie_{alg}^{\omega_Y} (Y)$ will be called rich
if (1) it generates  $\Lie_{alg}^{\omega_Y} (Y)$ as a Lie algebra
and (2)  $\VFA (Y)$ as a $\C [Y]$-module.
\end{definition}

\begin{remark}\label{pro2.remark}
(i) Note that (1) implies
the sets $\{ \delta_i (\C [Y]) \}$ generate the vector space $F_Y $.

(ii) For (2) it suffices to require that  the set of vector fields $\delta_1, \ldots , \delta_n, \ldots$
generates  the tangent space at each point of $Y$. This is a consequence of the next
simple fact  (e.g., see Exercise 5.8 in \cite{Ha}) which is essentially the Nakayama lemma.
\end{remark}

\begin{lemma}\label{vfi.no0}
 Let $A \subset B$ be a finitely
generated $\C [X ]$-module and its submodule. Suppose that for
every point $x \in X$ one has $A/M_x =B/M_x$ where $M_x$ is the
maximal ideal in $\C [X]$ associated with $x$. Then $A=B$
\end{lemma}

\begin{example}\label{frank.5}  Let $\sigma_1, \sigma_2, \ldots $ (resp.
 $\delta_1, \delta_2, \ldots $) be a a rich family on $X$ with respect to volume
 $\omega_X$ (resp.  on $Y$ with respect to volume
 $\omega_Y$). Denote their natural lifts to $X \times Y$ by the same
 symbols. Consider  the set $\cS$ of   ``horizontal" and ``vertical"
 fields of form $f\sigma_i$ and $g\delta_j$ where $f$ (resp. $g$) is a lift of a function
 on $Y$ (resp. $X$) to $X\times Y$.
 It follows from  the explicit construction in the proof of Proposition \ref{pro1} that $\cS$
 generates the Lie algebra $\Lie_{\omega} (X\times Y )$ for $\omega =\omega_X\times \omega_Y$,
 i.e. it satisfies condition (1) of Definition \ref{pro2}.
 Remark \ref{pro2.remark} (2) implies that condition (2) also holds and, therefore,
 $\cS$ is a rich family on $X \times Y$.

 In particular, consider a torus $\T = (\C^*)^n$ with coordinates $z_1, \ldots , z_n$.
 One can see that the vector field
 $\nu_j= z_j \partial /\partial z_j$ is a rich family on the $j$-th factor with respect
 to the invariant volume on $\C^*$. Thus fields of form $f_j\nu_j \,  (j=1, \ldots ,n)$ with
 $f_j$ being independent of $z_j$ generate
 a rich family on $\T$ with respect to the invariant volume.
\end{example}

\begin{convention}\label{vfi.no1}
Furthermore, we suppose that vector fields
$\delta_1, \ldots ,\delta_n, \ldots  $ form a rich family $\cS$ in $ \Lie_{alg}^{\omega_Y} (Y)$ and there are vector fields $\delta_1', \ldots , \delta_l' , \ldots  \in\VF_{alg}^{\omega_Y} (W,p)$ such that up to nonzero constant factors the set of  their restrictions to any fiber of $p$ contains $\delta_1, \ldots ,\delta_n, \ldots $ under some isomorphism between this fiber and $Y$.\footnote{If $p: W\to X$ is  a Zariski locally trivial fibration this Convention is automatically true. }
\end{convention}

\begin{lemma}\label{ffi.le2} \hfill

{\rm (1)} A function $g \in \C [W]$ is contained in $F(W,p)$ if and
only if its restriction to each general fiber of $p$ belongs to
$F_Y$. Furthermore, if $Y$ is fine then $\C [W] \simeq F(W,p) \oplus \C [X]$\footnote
{In fact for (1) one needs only that  the sets $\{ \delta_i (\C [Y]) \}$ generate the vector space $F_Y $
with $\delta_i$ running over the family $\cS$ from Convention \ref{vfi.no1} .}.

{\rm (2)} Suppose that
$Y$ has the algebraic volume density property. \break
Then
$\VF_{alg}^{\omega_Y}(W,p) = \Lie^{\omega_Y}(W,p)$.
\end{lemma}

\begin{proof}
There exists a cover $X = \bigcup_i X_i$ such that for each
$i$ one can find an \'etale surjective morphism $X_i' \to X_i$ for
which the variety $W_i': = W_i\times_{X_i}X_i'$ is naturally
isomorphic to $X_i' \times Y$ where $W_i =p^{-1}(X_i)$. Lifting
functions on $W$ to some $W_{i_0}'$ which is the direct product we
can introduce filtration $\C [W] = \bigcup_{i\geq 0} G_i$ as we did in the Claim in Proposition
\ref{pro1} (i.e. take a closed embedding $Y \hookrightarrow \C^m$
and consider the minimal degrees of extensions of function on
$W_{i_0}'$ to $X_{i_0}' \times \C^m$ with respect to the second
factor).

Consider now the set $S$ of functions  $g \in \C [W]$ such that  for a general $x \in X$ the restriction
$g|_{p^{-1}(x)}$ is in $F_Y$. Since the degree (generated by the
embedding $Y \hookrightarrow \C^m$) of the restriction of $g$ to
each (not necessarily general) fiber $p^{-1}(x_0) \simeq Y$ is
bounded by the same constant we see that $g|_{p^{-1}(x_0)}$
belongs to $F_Y$ (because the finite-dimensional subspace of $F_Y$
that consists its elements, whose degrees are bounded by this
constant, is closed).  Set $S_i=S\cap  G_i$ and $F_i=F(W,p)\cap G_i$.
If suffices to show that $S_i=F_i$ for every $i$.  Note that both $S_i$ and $F_i$ are finitely
generated $\C[X]$-modules and because of the existence of fields $\delta_1', \ldots , \delta_l',
\ldots $ in Convention \ref{vfi.no1}
 we see that for every $x \in X$ there is an equality $S_i/M_x =F_i/M_x$ where $M_x$ is the maximal ideal in $\C[X]$ associated with $x$. Hence the first statment  of (1) follows  from Lemma \ref{vfi.no0}.

For the second statement  of (1) note that $G_i/M_x =(F_i \oplus \C[X])/M_x$ since $Y$ is fine
and $\C [X]/M_x =\C$. Thus another application of Lemma \ref{vfi.no0} implies the desired conclusion.

The filtration of functions that we introduced, yields a filtrations of vector fields $\VF_{alg}^{\omega_Y} (W,p) = \bigcup B_i$ and $\Lie^{\omega_Y} (W,p)=
\bigcup L_i $ as where $B_i$ and $L_i$ are again finitely
generated $\C [X ]$-modules. Because $\cS$ in Convention \ref{vfi.no1}
is rich and $Y$ has the algebraic volume density property we have $B_i/M_x=L_i/M_x$. Thus  by Lemma \ref{vfi.no0} $B_i=L_i$
which yields (2).
\end{proof}

\begin{definition}\label{vfi.de1} (1) Suppose  that  Convention \ref{vfi.no1} holds,
$\omega_X$ is a volume form on $X$,
 $X = \bigcup X_i$, $ X_i'$, $W_i$, and
$W_i'$ are as in the proof of Lemma \ref{ffi.le2},  $\varphi_i : W_i'\to
X_i'\times Y$ is the natural isomorphism and $\omega$ is an
algebraic volume form on $W$ such that up to a constant factor
$ \varphi_i^*\omega$ coincides with $
((\omega_X \times \omega_Y)|_{X_i \times Y})$ for each $i$. Then
we call $p : W \to X$ a volume fibration (with respect to the
volume forms $\omega_X$, $\omega_Y$, and $\omega$).

(2) We call a derivation $\sigma' \in \ \Lie^{\omega} (W)$ a lift
of a derivation $\sigma
\in \Lie^{\omega_X} (X)$ if  for every $w \in W$ and $x=p(w)$ one
has $p_* (\sigma'(w)) =\sigma (x)$. (Note that  the Lie bracket of
two lifts is a lift.) We say that
$\sigma'$ is $p$-compatible if for any $\delta' \in \VF_{alg}^{\omega}
(W,p)$ we have $[\sigma' , \delta' ]\in  \VF_{alg}^{\omega} (W,p)$ and
the span ${\rm Span \,} (\Ker \sigma' \cdot \Ker \delta')$
coincides with $\C [W]$.

(3) We call a volume fibration $p : W \to X$ refined if the fiber $Y$ is a fine manifold with
 algebraic volume density property, and a rich family of fields on $X$ has $p$-compatible lifts.

We shall see later (Lemma \ref{refine.1}) that for a reductive group $G$ and its Levi semi-simple subgroup $L$
the natural morphism $G \to G/L$ can be viewed as a refined
 volume fibration with respect to appropriate volume forms.
\end{definition}

Since any algebraic vector field tangent to fiber of $p : W \to X$
has zero $\omega$-divergence if and only if its restriction on
each fiber has zero $\omega_Y$-divergence we have the following consequence
of Lemma \ref{ffi.le2} (2).

\begin{corollary} \label{vfi.co1} Let $p : W \to X$ be a volume
fibration whose fiber has the algebraic volume density property
 and $\VF_{alg}^{\omega}(W,p)$ (resp. $\Lie^{\omega}(W,p)$) be
the space of zero $\omega$-divergence algebraic vector fields
tangent to the fibers of $p$ (resp. the Lie algebra generated by
completely integrable algebraic vector fields tangent to the
fibers of $p$ and of zero $\omega$-divergence). Then
$\VF_{alg}^{\omega}(W,p) = \Lie^{\omega}(W,p)$.

\end{corollary}

The next fact will not be used further but it is interesting by itself.

\begin{proposition}\label{ffi.pro1} Let $p: W \to X$ be a refined volume fibration with
a fine fiber $Y$ and fine base $X$. Then $W$ is
fine.

\end{proposition}

\begin{proof} Indeed, the existence of lifts in Definition \ref{vfi.de1}  of a rich family
makes $F_X$ and, therefore, $F(W,p) \oplus
F_X$ a natural subspace of $F_W$. It remains to note that $\C [W]
= F(W,p) \oplus F_X \oplus \C$ by the assumption and by Lemma
\ref{ffi.le2}.

\end{proof}

%

\begin{proposition}\label{vfi.pro2} Let $p : W \to X$ be a refined volume fibration and
$\Theta$ be the
set of $p$-compatible lifts of a rich family from $\Lie_{alg}^{\omega_X}
(X)$ . Consider the space $L$ generated by $\Lie^{\omega}(W,p)$
and vector fields of form $\nu: =[f\sigma' , \delta' ]$
where $\sigma' \in \Theta  , \delta' \in \Lie^{\omega}(W,p)$, and
$f \in \Ker \sigma'$. Suppose that $T = p^* (TX)$ is the pull-back
of the tangent bundle $TX$ to $W$, $\varrho : TW \to T$ is the
natural projection, and $\cL = \varrho (L) $. Then

{\rm (1)} $\cL$ is a $\C [X]$-module;

{\rm (2)} $\cL$ consists of all finite sums $\sum_{\sigma' \in
\Theta} h_{\sigma'} \varrho (\sigma' )$ where $h_{\sigma'} \in
F(W,p)$.
\end{proposition}

\begin{proof}
The space $Lie^{\omega} (W,p)$ is, of course, a $\C [X]$-module.
Thus it suffices to consider fields like $\nu =[\sigma' ,\delta']$
only. Since $\sigma'$ is a lift of $\sigma \in \Lie^{\omega_X}
(X)$ we see that $\sigma' (\C [X] ) \subset \C[X]$ where we treat
$\C [X]$ in this formula as a subring of $\C [W]$. Then for every
$\alpha \in \C [X]$ we have $\alpha \nu = [\sigma', \alpha \delta'
] - \sigma' (\alpha)\delta'$ which implies (1).

Let $f_1 \in \Ker \sigma'$, $f_2 \in Ker \delta'$, and $h=\delta'
( f_1f_2)$. Then by the $p$-compatibility assumption $[f_1\sigma'
, f_2\delta' ] =h \sigma' + a$ with $a \in \VFA (W,p)$ and, furthermore, the span
of functions like $h$ coincides with $\delta' (\C [W])$. Thus
$\cL$ contains $F(W,p)\varrho (\sigma')$ which is (2).

\end{proof}

\begin{theorem}\label{ffi.th1}
Let $p : W \to X$  be a refined volume fibration whose fiber and base
have the algebraic volume density property.
Then $W$ has the algebraic volume density property.

\end{theorem}

\begin{proof} Suppose that $\delta_i'$ is as in
Convention \ref{vfi.no1}, $\{ \sigma_i \}$ is a rich family on $X$,
and  $\sigma_i' \in \Lie^{\omega} (W)$ is a $p$-compatible
lift of $\sigma_i$.  Let $\kappa \in \VF_{alg}^{\omega}
(W)$. Then $\kappa =\sum_{i} h_i \sigma_i' + \theta$ where $h_i
\in \C[W]$ and $\theta \in \VFA (W,p)$. By Proposition
\ref{vfi.pro2} and Lemma \ref{ffi.le2}(1)  adding an element of $L$
to $\kappa$ one can suppose that each $h_i \in \C[X]$. Since
$\theta$ is a $\C [W]$-combination of  $\delta_1', \ldots , \delta_l', \ldots  $
and $\diver_{\omega}( f \delta_i') = \delta_i' (f)$,
we see that $\diver_{\omega} \theta \in F(W,p)$. On the other hand
$\diver_{\omega}\sum_{i} h_i \sigma_i'=\sum_i \sigma_i (h_i)=
\diver_{\omega_X}\sum_{i} h_i \sigma_i \in \C [X]$. Hence
$\diver_{\omega_X}\sum_{i} h_i \sigma_i=0$ and $\theta \in
\VF_{alg}^{\omega} (W,p)$. By assumption $ \sum_{i} h_i \sigma_i
\in Lie_{alg}^{\omega_X} (X)$. In combination with the existence of
lifts for $\sigma_i$ and Corollary \ref{vfi.co1}
this implies the desired conclusion.

\end{proof}

\section{Volume forms on homogeneous spaces}

\begin{definition}\label{volume.1} We say that
an affine algebraic variety $X$ is (weakly) rationally connected
if for any (resp. general) points $x,y \in X$ there are a sequence
of points $x_0=x,x_1,x_2, \ldots , x_{n}=y$ and a sequence of
polynomial curves $C_1, \ldots ,C_n $ in $X$ such that
$x_{i-1},x_i \in C_i$.
\end{definition}

\begin{remark}\label{volume.2} Since finite morphisms transforms polynomial
curves into polynomial curves we have the following: if $X$ is an
affine (weakly) rationally connected variety and $f : X \to Y$ is
a finite morphism then $Y$ is also an affine (weakly) rationally
connected variety.

\end{remark}

\begin{example}\label{volume.3} It is easy to see that $SL_2(\C )$ is affine
rationally connected. (Indeed, $SL_2(\C )$ can be presented as an
algebraic locally trivial $\C$-fibration over $\C^2$ without the
origin $o$. Over any line in $\C^2$ that does not contain $o$ this
fibration is trivial and, therefore, admits sections which implies
the desired conclusion.)  Hence any semi-simple group is
rationally connected since its simply connected covering is
generated by $SL_2(\C )$-subgroups.

\end{example}

\begin{proposition}\label{volume.4}  Let $X$ be an affine manifold
and $\omega, \omega_1$ be algebraic
volume forms on $X$.

 {\rm (1)} If $X$ is weakly rationally
connected  then $\omega =c\omega_1$ for nonzero constant $c$.

{\rm (2)} If $G$ is a rationally connected linear algebraic group (say, unipotent or semi-simple)
acting on $X$ then $\omega \circ \Phi_g = \omega$ for the action
$\Phi_g : X \to X$ of any element $g \in G$.

\end{proposition}

\begin{proof} (1) Note that $\omega =h\omega_1$ where $h$ is an
invertible regular function on $X$. Let $x,y,x_i,C_i$ be as in
Definition \ref{volume.1}. By the fundamental theorem of algebra
$h$ must be constant on each $C_i$. Hence $h(x) =h(y)$ which
implies the first statement.

(2) Let $e$ be the identity in $G$. Then
we have a sequence $g_0=e, g_1, g_2, \ldots , g_m=g$ in $g$ such that
for any $i \geq 1$ there is a polynomial curve $C_i$ in $G$ joining $g_{i-1}$ and $g_i$.
Again for every $a \in C_i$ we have $\omega \circ \Phi_a =h(a) \omega$ where $h$ is
a nonvanishing regular function on $C_i$, i.e. a constant. This implies
 $\omega \circ \Phi_g =  \omega \circ \Phi_e = \omega$ which concludes the proof.
\end{proof}

For a Lie group $G$ one can construct a left-invariant (resp.
right-invariant) algebraic volume form by spreading the volume
element at identity by left (resp. right) multiplication (one of
these forms can be transformed into the other by the automorphism
$\varphi : G \to G$ given by $\varphi (g) =g^{-1}$).
Proposition \ref{volume.4} yields now the following well-known facts.

\begin{corollary}\label{volume.5} For a semi-simple Lie group $G$
its left-invariant volume form is automatically also right-invariant.
\end{corollary}

\begin{remark}\label{volume.5a} Since up to a finite covering
any reductive group $G$ is a
product of a torus and a semi-simple group we see that
the left-invariant volume form on this group is also right-invariant.

\end{remark}

\begin{proposition}\label{volume.6} Let $W$ be a linear algebraic group
group, $Y$ be its rationally connected subgroup,
and $X =W/Y$ be the
homogeneous space of left cosets. Then there exists  an algebraic
volume form $\omega_X$ on $X$ invariant under the action of $W$
generated by left multiplication.

\end{proposition}

\begin{proof}
Consider a left-invariant volume
form $\omega$ on $W$ and left-invariant vector fields $\nu_1, \ldots , \nu_m$ on
the coset $eY\simeq Y$, where $e$ is the identity of $W$,  so that they
generate basis of the tangent space at any point of this coset.
Extend these vector fields to $W$ using left
multiplication. Since $eY$ is a fiber of the natural projection $p : W \to X$ and left
multiplication preserves the fibers of $p$ we see that the
extended fields are tangent to all fibers of $p$. Consider the left-invariant form $\omega_X=
\iota_{\nu_1} \circ \ldots \circ \iota_{\nu_m} (\omega)$ (where $\iota_{\nu_i}$ is the
inner product of vector fields and differential forms). By construction it can be viewed
as a non-vanishing form on vectors from the pull-back of the tangent bundle $TX$ to $W$.
To see that it is actually a volume form on $X$ we have to show that it is invariant under right multiplication by
any  element $y \in Y$.  Such multiplication generates an automorphism of
$TW$ that sends vectors tangent (and, therefore,
transversal) to fibers of $p$ to similar
vectors. Hence it transforms $\omega_X$ into $f_y\omega_X$
where $y \to f_y$ is an algebraic homomorphism from $Y$ into the group of non-vanishing regular functions of $W$.
Since the rationally connected group $Y$ has no nontrivial
algebraic homomorphisms into $\C^*$ we have $f_y \equiv 1$ which yields the desired conclusion.

\end{proof}

By Mostow's theorem \cite{Mo} a linear algebraic group $W$ contains a Levi
reductive subgroup $X$ such that as an affine algebraic variety $W$ is
isomorphic to $X\times Y$ where $Y$ is the unipotent radical of $W$.
More precisely, each element $w \in W$ can
be uniquely presented as $w=g \cdot r$ where $g \in X$ and $r \in Y$. This presentation
allows us to choose this isomorphism $W \to X\times Y$ uniquely.

\begin{corollary}\label{refine.2} For the isomorphism  $W \to X\times Y$ as before
the left invariant volume form $\omega$ on $W$ coincides with $\omega_X \times \omega_Y$
where $\omega_X$ is a left-invariant volume form on $X$ and $\omega_Y$ is
an invariant form on $Y$.

\end{corollary}

\begin{proof} Note that $\omega$ is invariant by left multiplications (in particular by elements of $X$)
and also by right multiplication by elements of $Y$ (see, Lemma \ref{volume.4} (2)). This
determines $\omega$ uniquely up to a constant factor. Similarly, by
construction $\omega_X \times \omega_Y$ is invariant by left multiplications by elements of $X$
and by right multiplication by elements of $Y$.

\end{proof}

\begin{example}\label{refine.3} Consider the group $W$ of affine automorphisms
$z \to az +b$ of the complex line $\C$ with coordinate $z$. Then $Y\simeq \C_+$ is the group of
translations   $z \to z+b$ and we can choose $X\simeq \C^*$ so that its elements are automorphisms
of form $z \to az$. One can check that the left-invariant volume $\omega$ on $W$
coincides with ${\frac{{\rm d} a}{a^2}} \wedge {\rm d} b$ while $\omega_X ={\rm d} a/a$
and $\omega_Y ={\rm d} b$. The isomorphism $W \to X\times Y$ we were talking about
presents $az+b$ as a composition of $az$ and $z+b/a$.  Thus in this case
Corollary \ref{refine.2} boils down to the equality
$${\frac{{\rm d} a}{a^2}} \wedge {\rm d} b ={\frac{{\rm d} a}{a}} \wedge {\frac{{\rm d} b}{a}}.$$

\end{example}

\begin{proposition}\label{volume.7} Let $W$ be a linear algebraic group
group, $Y$ be its rationally connected subgroup,
and $X =W/Y$ be the
homogeneous space of left cosets. Suppose that $\omega$ is the left-invariant volume form on $W$, $\omega_Y$  the
invariant volume form on $Y$ and $\omega_X$ the volume form on the quotient constructed in Proposition \ref{volume.6}. Then the natural projection $p
: W \to X$ is a volume fibration with respect to the volume
forms $\omega , \omega_X$, and $\omega_Y$.

\end{proposition}

\begin{proof} Choose locally nilpotent derivations $\sigma_1',
\ldots , \sigma_k'$ and semi-simple derivations $\sigma_{k+1}'$, $\ldots , \sigma_n'$
on $W$ generated by the left multiplication of
$W$ by elements of its $\C_+$ and $\C^*$-subgroups and such that they
generate tangent space at each point of $W$. Since they commute
with morphism $p$ they yield locally nilpotent and semi-simple derivations
$\sigma_1, \ldots , \sigma_n$ on $X$ with the same property. Take
any point $x \in X$ and suppose that $\sigma_l, \ldots , \sigma_m$
generate the tangent space $T_xX$ (where $l \leq k \leq m$). Then we have the dominant
natural morphism $\psi : G \to X$ from the group $ G:=\C^{k-l+1}\times ( \C^*)^{m-k}$
given by the formula $\bar t =(t_l, \ldots , t_m) \to h_{t_l} \circ \cdots \circ h_{t_m} (x)$ where $\bar t \in G$ and
$h_{t_j}$ is the action of the element $t_j$ from the $\C_+$ or $\C^*$-group corresponding to the $j$-th factor in $G$.
This morphism
is \'etale at the identical element $ o = (0, \ldots , 0, 1, \ldots ,1)$ of $G$
and $\psi (o)=x$. The
restriction of $\psi$ to an open Zariski dense subvariety $Z$ of
$\C^{m-l+1}$ may be viewed as an \'etale neighborhood of $x \in X$.
Suppose that $\omega_Z$ (resp. $ \tilde{\sigma}_i$) is the lift of
the form $\omega_X$ (resp. vector field $\sigma_i$) to $Z$. By
Proposition \ref{volume.6} $\omega_Z$ is invariant under the local
phase flow generated by  $ \tilde{\sigma}_i$. Set $W'=W\times_X
Z$. Then by construction, $W'$
is naturally isomorphic to $Z \times Y$ and under this
isomorphism each field $\sigma_i'$ corresponds to the horizontal
lift of $\tilde{\sigma}_i$ to $Z \times Y$. Hence $\omega_Z \times
\omega_Y$ is invariant under the local phase flow generated by
this lift of $\tilde{\sigma}_i$. It is also invariant under right
multiplication by elements of $Y$ by Proposition \ref{volume.4} (2) and,
therefore, determined uniquely by its value at one point. But the
form $\omega$ is also invariant under the local phase flow
generated by ${\sigma}_i'$ and under right multiplication by
elements of $Y$ again by  Proposition \ref{volume.4} (2). Therefore, the preimage of
$\omega$ on $W$ coincides with $\omega_Z \times \omega_Y$ since one can choose
$\omega_Y$ so that both forms coincide at one point.

\end{proof}

We finish this section with the following example of a refined volume fibrations.

\begin{lemma}\label{refine.1} Let $W$ be a reductive group and $Y$ be its Levi
semi-simple subgroup. Then $p : W \to X: =W/Y$ is a refined volume fibration.

\end{lemma}

\begin{proof}  By Proposition \ref{volume.7} it is a volume fibration. Thus by
Definition \ref{vfi.de1} it remains to find a rich family of algebraic fields on $X$
that have $p$-compatible lifts. Let $T$ be the connected component of the center
of $G$. That is, $T$ is a torus $(\C^*)^n$ and $X \simeq T/(T\cap Y)$ is also a torus $(\C^*)^n$
since the group $T\cap Y$ is finite. Let us
start with the case when $T \cap Y$ is trivial, i.e. $X=T$. Then $W =X \times Y$
and $p$ is the projection to the first factor. In particular, any ``vertical" field
$\delta' \in \AVF (W,p)$ contains $\C [X ] \subset \C [W]$ in its kernel. Every vector field
 $\sigma \in \Lie_{\omega_X} (X)$
has  similarly a lift $\sigma' \in \Lie_{\omega} (W)$ such that this lift contains $\C [Y]$ in its kernel.
In particular, $\Span \Ker \sigma' \cdot \Ker \delta' = \C [W]$. Furthermore, any vertical
field $\delta'$ is of form $\sum_jf_j\delta_j$ where $f_j \in \C [W]$ and $\delta_j$ is
the natural  lift of vector field on $Y$ to $W$.
Since $[\sigma' , \delta_j ]=0$ we have $[\sigma', \delta'] \in \AVF (W,p)$
which shows that any rich family of vector fields on $X$
has the desired $p$-compatible lifts.

In the general case when $T\cap Y$ is not trivial we have a commutative
diagram
\[ \begin{array}{ccc}
T\times Y & \stackrel{{ \phi}}{\rightarrow} & W\\
\, \, \, \, \downarrow  q &   & \, \, \,  \downarrow p \\
T & \stackrel{{ \psi}}{\rightarrow} & X
\end{array}  \]
where the vertical arrows are unramified finite coverings.
Let $z_1, \ldots , z_n$ be natural
coordinates on $T\simeq( \C^*)^n$ and  $w_1, \ldots , w_n$ be natural
coordinates on $X\simeq (\C^*)^n$. Then up to constant factors we have
$w_j =\prod_{i=1}^nz^{k_{ij}}$.
By Example \ref{frank.5} a rich family on $X$ consists of vector fields of form
$\nu = f_j w_j \partial /\partial w_j$ where $f_j$ is a function on $X$ independent of $w_j$.
Note that $w_j \partial /\partial w_j = \sum_{i=1}^n k_{ij}z_i\partial /\partial z_i$.
Since $z_i\partial /\partial z_i$
is associated with multiplication by elements of a  $\C^*$-subgroup of $T$
it may be viewed as a field on $X$ and we can find its lift to $W$. Thus the
fields $w_j\partial /\partial w_j$ and also $\nu$ have  lifts to $W$ and we need to check
that they are $p$-compatible.

Let $\sigma'$ be the lift of one of these fields and $\sigma''$ be its
preimage on $T\times Y$. Each vertical vector field $\delta'$
on $W$ (i.e. it is from the kernel of $p_*$) generates a vertical vector field $\delta''$ on
$T\times Y$ (i.e. it is from the kernel of $q_*$). As we showed  before
$[\sigma'',\delta'']\in \AVF (T\times Y, q)$ and therefore $[\sigma' ,\delta']\in \AVF (W,p)$.
Furthermore, since $\Span \Ker \sigma'' \cdot \Ker \delta'' = \C [T \times Y]$
we have still the equality
$\Span \Ker \sigma' \cdot \Ker \delta' = \C [W]$ by virtue of Lemma \ref{app.6}. This
concludes the proof of $p$-compatibility and the Lemma.

\end{proof}

\section{Compatibility}

\begin{notation}\label{comp.1}
Let $G$ be a semi-simple Lie group, $S_0$ and $S$ be its $SL_2 $
or $PSL_2$-subgroups, and $p : G \to X :=G/S_0$ be the natural projection
into the set of left cosets. Suppose that $\delta$ is a completely
integrable algebraic vector field on $S_0$ generated  by right
multiplications. Then it generates $\delta' \in \LieAO (G,p)$. Let $H \simeq \C_+$ be
a subgroup of $S$. Left multiplication by elements of $H$
generate a locally nilpotent derivation $\sigma'$ on $G$. Note
that $[\sigma' ,\delta' ]=0$ (i.e. we have an $(S\times
S_0)$-action on $G$) and $\sigma'$ generates a locally nilpotent
derivation $\sigma$ on $X$ associated with the corresponding
$H$-action on $X$.
\end{notation}

We think of $S_0$ being fixed  and
our aim is to find "many" $S$ such  that $\sigma'$ is $p$-compatible for
 $S$, i.e. the vector space generated by $\Ker \sigma' \cdot \Ker
\delta'$ coincides with $\C [G]$. From now on we use the (seemingly overloaded) notation of strictly semi-compatibility
for pairs of vector fields (for Definition see the Appendix) since it was introduced in the work of \textsc{Donzelli, Dvorsky} and the first author \cite{DoDvKa} and we like to stick to this earlier introduced notation. We apologize for any inconvenience to the reader.

\begin{lemma}\label{comp.2} Suppose $g_0\in G$ and $S \cap g_0S_0g_0^{-1} = \Gamma$.
Then the isotropy group of the point $g_0S_0 \in X$ under the
$S$-action is $\Gamma$. In particular, if the $S$-orbit of
$g_0S_0$ is closed then $\Gamma$ is reductive by the Matsushima
theorem.
\end{lemma}

\begin{proof} The coset $g_0S_0$ is fixed
under the action of $s \in S$ if and only if $sg_0 S_0 \subset
g_0S_0$ which implies that $g_0^{-1}sg_0\in S_0$ and we have the
desired conclusion.

\end{proof}

\begin{proposition}\label{comp.5} Let $\Gamma_g = S \cap gS_0g^{-1}$ be finite
for every $g \in G$. Then $\sigma'$ is $p$-compatible.
\end{proposition}

\begin{proof} Consider the quotient morphism $r: G \to Z :=G//(S \times
S_0)$. Since $\Gamma_g$ is always finite all orbits are
equidimensional and, therefore, closed (indeed, for a reductive
group $S\times S_0$ the closure of a non-closed orbit must contain
a closed orbit, automatically of smaller dimension, which is impossible because all orbits are of the same dimension).
By Luna's slice theorem  for
every point $z \in Z$ there exists a Zariski neigborhood $U\subset Z$, a
$\Gamma_g$-invariant slice $V \subset G$ through a point of
$r^{-1}(z)$ such that $r|_V : V \to U$ is a surjective
quasi-finite morphism, and a surjective \'etale morphism $W \to
r^{-1} (U)$ where $W = V\times_{\Gamma_g} (S \times S_0)$. In
particular, we have a natural surjective quasi-finite morphism
$W'':=V\times (S \times S_0) \to r^{-1}(U)$. Clearly, the
algebraic vector fields $\sigma''$ and $\delta''$ on $W''$ induced
by $\sigma'$ and $\delta'$ are strictly semi-compatible,
i.e. the span of $\Ker \sigma'' \cdot  \Ker \delta''$ coincides
with $\C [W'']$. Note also that for any $\C_+\simeq H < S$ the
quotient $G//H$ is smooth and the quotient morphism $G \to G//H$
is a holomorphic $\C$-fibration over its image. By Lemmas
\ref{app.6} and \ref{app.8} in Appendix the restrictions of
$\sigma'$ and $\delta'$ to $r^{-1}(U)$ are also strictly
semi-compatible. Thus there is a cover $Z = \bigcup U_i$ such that each
$U_i$ is of form $U_i =Z \setminus g_i^{-1}(0)$ with $g_i \in \C
[Z]$ and the restrictions of of $\sigma'$ and $\delta'$ are
strictly semi-compatible on each $W_i=r^{-1}(U_i)$. For any function
$h \in \C[G]$ its restriction $h\vert_{W_i}$ is contained in $\Ker \sigma' \vert_{W_i}\cap \Ker \delta'\vert_{W_i}$. Since for any function $\varphi \in \C[W_i]$ there exists $m>0$ such that $\varphi g_i^m \in \C [G]$ and
since $g_i\in \Ker \sigma' \cap \Ker \delta'$, for an appropriate
$m$ the function $hg_i^m$ belongs to the span of $\Ker \sigma'
\cdot \Ker \delta'$. Now the desired conclusion follows from the
standard application of the Nullstellensatz.

\end{proof}

\begin{remark}\label{comp.6} Lemmas \ref{app.6} and \ref{app.8} are
to a great extend repetitions of Lemmas 3.6 and Lemma 3.7 in
\cite{KaKu1}. Therefore, we put them in Appendix.
\end{remark}

\begin{lemma}\label{comp.8} Let $G,S_0,X$, and $S$  be as in Lemma
\ref{comp.2} and $\Gamma_g =S \cap gS_0g^{-1}$ where $g \in G$.
Suppose that
$\Gamma_g$ does not contain a torus $\C^*$ for every $g \in G$.
Then every $\Gamma_g$ is finite.

\end{lemma}

\begin{proof} Assume
that $\Gamma_{g_0}$ is not finite for some $g_0 \in G$. Then
$\Gamma_{g_0}$ cannot be reductive (without a torus) and the
$S$-orbit $O$ of $g_0S_0 \in X$ is not closed by the second
statement of Lemma \ref{comp.2}. Furthermore, since any
two-dimensional subgroup of $SL_2( \C )$ contains $\C^*$ we see
that $\Gamma_{g_0}$ is one-dimensional, i.e. $O$ is
two-dimensional. Since $S$ is reductive the closure of $O$ must
contain a closed orbit $O_1$ of some point $g_1S_0 \in X$. Thus
$\dim O_1 \leq 1$ and $\dim \Gamma_{g_1} \geq 2$. But in this case
as we mentioned $\Gamma_{g_1}$ contains a torus which yields a contradiction.

\end{proof}

In order to find $S$ such that $\Gamma_g=g^{-1}S_0 g \cap S$ does not
 contain a torus for every $g \in G$ we need to remind the notion of a principal $SL_2$ or
 $PSL_2$-subgroup
 of a semi-simple group  $G$
 (resp. principal $\sgoth \lgoth_2$-subalgebra in the Lie algebra $\ggoth$ of $G$) from \cite{Bo}.  Recall that a semi-simple element $h$ of $\ggoth$ is called regular if the dimension of  its centralizer is equal to the rank of $\ggoth$ (more precisely, this centralizer coincides with a Cartan subalgebra
 $\hgoth$ of $\ggoth$).  An $\sgoth \lgoth_2$ subalgebra $\sgoth$ of $\ggoth$ is called principal if it contains a
 regular semi-simple element $h$ such that every eigenvalue of its adjoint operator is an even integer.
The $SL_2$ (or $PSL_2$) subgroup generated by such subalgebra is also called principal. For instance, in $SL_n$ up to conjugation every regular element is a diagonal matrix with distinct eigenvalues and any principal $SL_2$-subgroup acts irreducibly on the natural $n$-space. Any two principal $SL_2$-subgroups are conjugated in $G$ and any $SL_2$-subgroup corresponding to a root is not principal (unless $\ggoth =\sgoth \lgoth_2$) since its semi-simple elements are not regular.

\begin{lemma} \label{comp.10} If $S$ is a principal  $SL_2$ (resp. $PSL_2$) subgroup of a semi-simple group $G$ and $S_0$ be any subgroup of $G$ that does not contain regular semi-simple elements. Then $\Gamma_g =g^{-1}S_0g \cap S$ is finite for every $g \in G$.

\end{lemma}

\begin{proof}  Note that $\Gamma_g$ cannot contain a torus since otherwise $S_0$ contains a regular semi-simple element.  Lemma \ref{comp.8} implies now the desired conclusion.

\end{proof}

\begin{proposition}\label{comp.20} Let $G$ be a
semi-simple Lie group different from $SL_2(\C )$ or $PSL_2(\C )$.
Suppose that $S_0,Z=G/S_0$, $p: G\to Z$, and $\sigma'$
are is in Notation \ref{comp.1}. Let $S_0$ correspond to a root in the
Dynkin diagram. Then $\sigma'$ can be chosen that
it is $p$-compatible (for any $S_0$ corresponding to a root in the Dynkin diagram!!).
 Furthermore, there are enough of these $p$-compatible completely
integrable algebraic vector fields $\sigma'$, so that the Lie
algebra $L$ generated by them generates $\VFA (Z)$ as a $\C [Z]$-module.

\end{proposition}

\begin{proof} Let an $SL_2$ (or $PSL_2$) subgroup $S_{0}$ correspond to a
root  and $S$ be a principal $SL_2$ (or $PSL_2$) subgroup. By Proposition \ref{comp.5}
and Lemma  \ref{comp.10}    $\sigma'$ is $p$-compatible and we are left with the second
statement.  Suppose that $X,Y,H$ is a standard triple in the $\sgoth \lgoth_2$-subalgebra $\sgoth$ of $S$, i.e. $[X,Y]=H$, $[H,X]=2X$,  $[H,Y]=-2Y$ \footnote{It is unfortunate,
but we have to use the classical notation $X,Y,H$ for a standard triple of an $\sgoth \lgoth_2$-algebra while in the  rest of the text these symbols denote affine algebraic varieties and groups.}.
In particular, the locally nilpotent vector fields generated by $X$ and $Y$ are of form $\sigma'$ and they are $p$-compatible. Suppose that the
centralizer of $H$  is
the Cartan subalgebra  $\hgoth$ associated with the choice of  root system and $X_0,Y_0,H_0$ is an $\sgoth \lgoth_2$-triple corresponding to one of the roots. Conjugate $S$ by  $x_0=e^{\varepsilon X_0}$ where $\varepsilon$ is a small parameter. Up to terms of order 2 element $H$ goes to $H+\varepsilon [H,X_0]$ after such conjugation, i.e. $[H,X_0] $ belongs (up to second order) to  the Lie algebra generated by $X,Y$, and the nilpotent elements of the Lie algebra of principal $SL_2$-subgroup $x_0^{-1}Sx_0$. Since each $X_0$ is an eigenvector of the adjoint action of $H$ we have $[H,X_0]=aX_0$. Furthermore, $a \ne 0$ since otherwise $X_0$ belongs to the centralizer $\hgoth$ of the regular element $H$. Thus $X_0$ and similarly $Y_0$ are (up o second order) in the Lie algebra $L$ generated by fields of form $\sigma'$. The same is true for $H_0=[X_0,Y_0]$.  Thus values of $L$ at any point $z \in Z$ generate the tangent space $T_zZ$
which implies $L$ generates $\VFA (Z)$ as a $\C [Z]$-module.
\end{proof}

\section{Main Theorem}

\begin{notation}\label{main.1} In this section $G$ is a
semi-simple Lie group.
By $S_i$  we denote an $SL_2$ or $PSL_2$-subgroup of $G$ (for each index $i\geq 0$)
and by  $p_i : G \to X_i =G/S_i$ the natural projection.
By abusing notation
we treat $\C [ X_i]$ as the subring $p_i^* (\C [X_i])$ in $\C [G]$ . Note that Lemma
\ref{ffi.le2} implies that $\C [G ] \simeq F(G,p_i) \oplus \C [X_i ]$ and denote by $\pr_i  : \C [G ] \to \C [X_i]$ the natural projection.
For any semi-simple
complex Lie group $B$ denote by $B^{\R}$ its maximal compact
subgroup whose complexification coincides with $B$ (it is unique up to conjugation).  Let
$K_i=S_i^{\R}$. Define a linear operator $\av_i : \C [G]
\to \C [G]$ by
$$\av_i (f)= \int_{K _i} f(wk) \, {\rm d} \, \mu_{K_i} (k)$$ for any function $f \in \C [G]$ where $\mu_{K_i}
(k)$ is the bi-invariant normalized Haar measure on $K_i$.

\end{notation}

\begin{lemma}\label{main.2} In Notation \ref{main.1} we have

{\rm (i)} the right multiplication by an element $k \in K_i$
generates a map $\Psi : \C [G ] \to \C [G]$ (given by $f(w) \to
f(wk)$) whose restriction to $F(G,p_i)$ is an isomorphism;

{\rm (ii)} $\Ker \av_i =F(G,p_i)$, i.e. $\av_i =\pr_i$ and
$f -\av_i (f) \in F(G,p_i)$ for every $f \in \C [G]$.
\end{lemma}

\begin{proof} The right multiplication transforms every fiber
$Y:=p_i^{-1}(x)$ into itself and each completely integrable
algebraic vector field on it into a similar field. Hence for every
$f \in F(G,p_i)$ we have $\Psi (f)|_{p_i^{-1}(x)} \in F_Y$. Now (i)
follows from Lemma \ref{ffi.le2}. Thus operator $\av_i$ respects the
direct sum $\C [G ] \simeq F(G,p_i) \oplus \C [X_i ]$ and sends
$\C [G]$ onto $\C [X _i]$ so that its restriction to $\C [X_i
]$ is identical map. This implies (ii).

\end{proof}

\begin{lemma}\label{main.4} Let $S_0$ and $K_0$ be as before and let $L = G^{\R}$
contain $K_0$.
Consider the natural inner product on $\C [G] $ given by
$$h_1 \cdot h_2 =\int_{l \in L}h_1 (l) {\bar h_2} (l) {\rm d}
\mu_{L} (l)$$ where $\mu_{L}$ is the bi-invariant measure on
${L}$.  Then $\C [X_0]$ is the orthogonal complement of $F(G,p_0)$.
\end{lemma}

\begin{proof}
Consider $h_1 \in \C [G])$ and $h_2 \in \C [X_0]$.
Show that ${\rm av_0} (h_1) \cdot h_2=h_1 \cdot h_2$. We have  $$I:={\rm av_0} (h_1)\cdot h_2=
\int_{L}\int_{K_0}  h_1(lk_0)
{\bar h_2} (l)  {\rm d} \mu_{K_0} (k_0) {\rm d}
\mu_{L} (l).$$ By Fubini's theorem $$I= \int_{K_0}
\int_{L} h_1(lk_0) {\bar h_2} (l){\rm d} \mu_{L} (l)
 {\rm d} \mu_{K_0} (k_0) .$$  Set
$l'=lk_0$. Then $h_1(lk_0)=h_1(l')$ and
$h_2(l)=h_2(l'k_0^{-1})=h_2(l')$ since $h_2$ is right $K_0$-invariant. Using the fact
that measures are invariant we see that $I$ coincides with
$$ \int_{K_0}\int_{L} h_1(l') {\bar h_2}
(l') {\rm d} \mu_{L} (l') {\rm d}
\mu_{K_0} (k_0)= $$
$$\int_{L}h_1(l') {\bar h_2} (l')){\rm d} \mu_{L} (l')$$ where the last
equality holds since measure $\mu_{K_0}$ is normalized.  Thus ${\rm av} (h_1) \cdot h_2= h_1 \cdot h_2$.
Now the desired conclusion follows from Lemma \ref{main.2} and the fact that $\C [G] \simeq F(G,p_0) \oplus \C [ X_0]$.
\end{proof}

\begin{corollary}\label{main.6}  Let $S_0, \ldots , S_m$ be as in Notation \ref{main.1}
with each $K_i =S_i^{\R} \subset L$. Set  $F= \sum_{i=0}^m F(G, p_{i })$. Then the orthogonal
complement of $F$ in $\C [G]$ coincides with the subspace $V$ of functions  that are invariant
with respect to any $S_i$-action generated by the right multiplication. In particular, if this set
$S_0, \ldots , S_m$ contains all $SL_2$ or $PSL_2$-subgroups corresponding to positive (or even simple positive) roots for the Lie algebra of $G$ then this orthogonal complement  $V$ consists of constants only and $\C [G]\simeq F \oplus \C$.

\end{corollary}

\begin{proof} Indeed, treating $\C [X_i]$ as $p_i^*(\C [X_i])$ we see that by Lemma \ref{main.4} the orthogonal complement of $F$ in $\C [G]$ is $V=\bigcap_{i=0}^m \C [X_{i}]$ which is exactly the
space of functions invariant under each $S_i$-action.  For the second statement note that if the sequence $\{ S_i \}$ of subgroups generate the whole group $G$ then these invariant functions must be constants.
\end{proof}

\begin{lemma}\label{main.62}  Let $\{ S_i \}_{i=0}^m$, $F$, and $V$ be as in Corollary \ref{main.6} and $f_0 \in \C [G] \setminus F$. Consider the smallest subspace $U\subset \C [G]$ that contains $f_0$ and such that for every $i$ and every $f\in U$ function $\av_i (f)$ is also in $U$. Then

{\rm (1) } $U$ is of some finite dimension $N$;

{\rm (2)} $\dim U \cap F = N-1$ and $\dim U \cap V =1$.

\end{lemma}

\begin{proof}  Consider a closed embedding $\rho : G \hookrightarrow \C^n$  such that the induced action of $G$ on $\C^n$ is linear. This yields a filtration on $\C [G]$ defined  by minimal degrees of polynomial extensions of regular functions on $G$ to $\C^n$.  Let $W_k$ be the subspace of $\C [G]$ that consists of functions of degree at most $k$ and $\Phi_l :  \C [G] \to \C [G]$ be the automorphism  given by $f(w) \to f(wl)$ for $l \in L$.  Since the $G$-action on $\C^n$ is linear each automorphism  $\Phi_l$ sends $W_k$ into itself. Hence the definition of $\av_i$ implies that  $\av_i (W_k)\subset W_k$. Thus $U\subset W_k$ as soon as $f_0 \in W_k$ which yields (1).

 Denote the orthogonal projection onto $V$ by $\pr :  \C[G]  \to V$ and let $f_0'=\pr (f_0)$. Since $f_0 \notin F$ we have $f_0'\ne 0$.  Let $P$ be the hyperplane in $\C [G]$ that consists of vectors of form $f_0'+P_0$ where $P_0$ is the hyperspace orthogonal to $f_0'$.  In particular $P$ contains $f_0$. Since $f_0' \in \C [X_i]$ for every $i$ we see that $P$ is orthogonal to each $\C [X_i]$. Recall that the operator $\av_i =\pr_i$ is just the orthogonal projection to $\C [X_i]$, i.e. $P$ is invariant with respect to these operators. In particular, if we set $f_J=\av_{j_1} \circ \cdots \circ \av_{j_s} (f_0)$  for a multi-index $J=(j_1, \ldots , j_s)$
with $j_t \in \{ 0, \ldots , m \}$  then $f_J \in P \cap U$.

 We want to show that for some sequence of such multi-indices  $f_J$ is convergent to a nonzero
 element of $V$  or, equivalently, $f_J'$ is convergent to an element of $V\cap P_0$ for $f_J'=f_J-f_0'$. Consider the subspace $U'$ generated by vectors of form $f_J'$. Let $U_i' = U' \cap \C [X_i ]$  and $I=(J,i)$, i.e. $f_I'=\av_i (f_J')$.   By construction the operator $\pr_i|_{U'}=\av_i|_{U'}$ is just the orthogonal projection to $U_i'$.  Hence if $f_J' \notin U_i'$ we have $||f_I'|| <||f_J'||$. Since $U'$ is finite-dimensional this implies that one can choose $\{ f_J \}$ convergent to an element $v \in U'$ and we can suppose that $v$ has the smallest possible norm. Then $\pr (v) =  v$ because of the last inequality. On the other hand
 $\pr  ( f_J' ) = \pr (f_J) - \pr (f_0') = \pr \circ \av_{j_1} \circ \cdots \circ \av_{j_s} \circ (f_0) - \pr (f_0) = f_0' - f_0' = 0$.
 Thus $v=0$. This shows $f_0' \in V\cap U$ and therefore  $\dim U \cap V \geq 1$.

 On the other hand $U$ contains a subspace $U_0$ generated by vectors of form $f_0 -f_J$. One can see that $U_0$ is of codimension 1 in $U$. Furthermore, $f_0-f_I=(f_0-f_J)+ (f_J-f_I)$. Note that $f_J-f_I=f_J - \av_i (f_J) \in F(G,p_i)\subset F$ by Lemma \ref{main.2}. Thus, using induction by the length of the multi-index $J$ one can show that $f_0 -f_J \in F$. That is, $\dim U \cap F \geq N-1$ which concludes the proof.

\end{proof}

\begin{proposition}\label{main.proposition}
 Any semi-simple group $G$ has the algebraic volume density with
respect to the invariant volume.

\end{proposition}

\begin{proof}
Choose $S_0,
\ldots , S_m$ as in Corollary \ref{main.6} and such that they
correspond to simple nodes in the
Dynkin diagram (it is possible since every semi-simple group $G$ has a compact real form, i.e.
we can suppose that $S_i^\R =K_i \subset L =G^\R$). Consider the natural projections $p_i : G \to X_i :=
G/S_i$ to the sets of left cosets.  Choose $p_i$-compatible
completely integrable algebraic vector fields $\sigma'$ as in
Proposition \ref{comp.20} and denote their collection by $\Theta$.
That is,
vector fields from $\Theta$ are of zero divergence, they commute with any
$\delta \in \VF_{alg}^{\omega} (G,p_i)$, and  they are independent from index $i$.
Furthermore, these fields from $\Theta$ can be viewed also as zero divergence
vector fields on $X_i$ that generate $\VFA (X_i)$
as a $\C [X_i]$-module.

 Let us fix an index $i$.
Any algebraic vector field on
$G$ is of form $$ \nu = \sum_{\theta \in \Theta} h_{\theta} \theta
+ \delta$$ where the sum contains only finite number of nonzero terms,
$h_{\theta} \in \C [G]$ and $\delta \in \VFA
(G,p_i)$. Since $SL_2(\C )$ and $PSL_2(\C )$ are fine, $\C [G] = F(G,p_i) \oplus \C [X_i]$
by  Lemma \ref{ffi.le2}.
Thus by virtue of Proposition \ref{vfi.pro2} adding fields from $\Lie_{alg}^{\omega} (G)$ to $\nu$ we get a field
$$ \nu_i = \sum_{\theta \in \Theta} h_{\theta}^i \theta
+ \delta_i$$ where $\delta_i \in \VFAO (G, p_i)$,
$h_{\theta}^i = \av_i (h_{\theta})$, and $\av_i=\pr_i$ is
as Notation \ref{main.1}.
That  is, $h_{\theta}^i \in \C [X_i]$.
Suppose that $\diver_{\omega} \nu =0$ and, therefore, $\diver_{\omega} (\nu_i )=0$.
Note that
$$\diver_{\omega}(\sum_{\theta \in \Theta} h_{\theta}^i \theta )
=\sum_{\theta \in \Theta} \theta (h_{\theta}^i) \in \C [X_i]$$ while
$\diver_{\omega} \delta_i \in F (G,p_i)$. Hence
$\diver_{\omega}(\sum_{\theta \in \Theta} h_{\theta}^i \theta
)=\diver_{\omega}(\delta_i)=0$.  Since $SL_2(\C )$ and $PSL_2(\C )$  have the
algebraic volume density property, $\delta_i
\in \Lie^{\omega} (G,p)$ by Corollary \ref{vfi.co1}.
Thus $\nu -\tilde \nu_i \in \Lie_{alg}^{\omega}(G)$
where
$$\tilde \nu_i = \sum_{\theta \in \Theta} h_{\theta}^i \theta .$$
In particular it suffices to show that $\tilde \nu_i
\in \Lie_{alg}^\omega (G)$ and, therefore, we can suppose that $\delta=0$ in
the original formula for $\nu$.  Repeating this procedure
we see that  $\nu -\tilde \nu_J \in \Lie_{alg}^{\omega}(G)$
where
$$\tilde \nu_J = \sum_{\theta \in \Theta} h_{\theta}^J \theta $$
for a multi-index $J=(j_1, \ldots , j_s)$
with $j_t \in \{ 0, \ldots , m \}$ and
 $h_{\theta}^J=\av_{j_1} \circ \cdots \circ \av_{j_s} (h_\theta )$.

By Corollary \ref{main.6} and Lemma \ref{main.62} the vector space
generated by  $h_{\theta}$ and functions of form $h_{\theta}^J$ is also
generated by constants and functions of form $h_{\theta} - h_{\theta}^J$.
Thus adding to  $ \nu = \sum_{\theta \in \Theta} h_{\theta} \theta$ vector fields of form
$\nu -\tilde \nu_\theta^J$ and $c \theta$ (where $c \in \C$) we can reduce the number of
nonzero terms in this sum. Hence $\nu \in \Lie_{alg}^\omega (G)$ which
implies the desired conclusion.

\end{proof}

\subsection{Proof of Theorem \ref{theorem}.}  Let us start with the case when
$G$ is reductive. Suppose that $Y$ is its Levi semi-simple subgroup.
Then by Proposition \ref{refine.1} $p : G\to X:= G/Y$ is a refined volume fibration
and by Theorem \ref{ffi.th1} $G$ has the algebraic volume density property.

Now consider an arbitrary linear algebraic group $G$ and let $Y$ be its
unipotent ideal and $X$ be a Levi reductive subgroup of $G$. By Corollary
\ref{refine.2} the Mostow isomorphism $G \to X\times Y$ makes the left invariant
volume $\omega$ on $G$ equal to $\omega_X\times \omega_Y$ where $\omega_X$
is left invariant on $X$ and $\omega_Y$ is invariant on $Y$. Now by
Proposition \ref{pro1} $G$ has the algebraic volume density property
with respect to $\omega$ which concludes the proof of our Main Theorem.

\begin{remark}  Theorem \ref{theorem} remains, of course, valid if instead
of the left invariant volume form we consider the right invariant one,
because the affine automorphism $G \to G, \, g\to g^{-1}$ transforms
the left invariant volume form into the right one while preserving the complete
integrability of the algebraic vector fields.

\end{remark}

\section{Appendix: Strictly semi-compatible
fields}

\begin{notation}\label{H_i} In this section
$H_i$ is isomorphic to $\C_+$ for $i=1,2$. We suppose also that
$X$ is a normal affine algebraic variety equipped with nontrivial
algebraic $H_i$-actions (in particular, each $H_i$ generates an
algebraic vector field $\delta_i$ on $X$). The categorical
quotients will be denoted $X_i=X//H_i$ and the quotient morphisms
by $\rho_i : X \to X_i$.
\end{notation}

\begin{definition}\label{app.2}
A pair $(\delta_1, \delta_2)$  of algebraic vector fields
(as in Notation \ref{H_i}) is called strictly semi-compatible  if the
span of $\Ker \delta_1 \cdot \Ker \delta_2$ coincides with $\C
[X]$.
\end{definition}



We shall need the following obvious geometric reformulation of
Definition.

\begin{proposition}\label{app.4}  Let $\delta_1$ and $\delta_2$
be as in Notation \ref{H_i}. Set $\rho = (\rho_1, \rho_2) : X \to
Y:= X_1 \times X_2$ and $Z$ equal to the closure of $\rho (X)$ in
$Y$. Then $\delta_1$ and $ \delta_2$ are strictly semi-compatible if
and only if $\rho : X \to Z$ is an isomorphism.

\end{proposition}

\begin{lemma}\label{app.6}
Let $X,H_i,X_i, \delta_i$, and $\rho_i$ be as in Notation
\ref{H_i} with $\delta_1$ and $\delta_2$ being strictly
semi-compatible and $[\delta_1 , \delta_2]=0$.  Set $\Gamma = H_1
\times H_2$. Let $X'$ be a normal affine algebraic variety
equipped with a non-degenerate $\Gamma$-action and $p : X \to X'$
be a finite $\Gamma$-equivariant morphism (for each $i=1,2$), i.e.
we have commutative diagrams
\[ \begin{array}{ccc}
X & \stackrel{\rho_i}{\rightarrow} & X_i\\ \, \, \, \downarrow p&&
\, \,  \, \, \downarrow q_i \\ X' &
\stackrel{\rho_i'}{\rightarrow} &
X_i'\\
\end{array} \]
\noindent where $\rho_i' : X' \to X_i'=X'//H_i$ is the quotient
morphism of the $H_i$-action on  $X'$. Suppose also that $\rho_1'$
makes $X'$ an  \'etale locally trivial $\C$-fibration over $\rho_1' (X')$, and $X_1, X_2$ are
affine\footnote{In all cases we apply this Lemma the $\C_+$-action
generated by $H_i$ extends to an algebraic $SL_2(\C )$-action and,
therefore, $X_i$ is affine automatically by the Hadziev theorem
\cite{Had}.}. Then ${\rm Span} (\C [X_1'] \cdot \C [X_2'])=\C
[X']$.

\end{lemma}

\begin{proof}
Since $p$ is finite, every $f \in \C [X_i] \subset \C [X]$ is a
root of a minimal monic polynomial with coefficients in $\C [X' ]$
that are constant on $H_i$-orbits (since otherwise $f$ is not
constant on these orbits). By the universal property these
coefficients are regular on $X_i'$, i.e. $f$ is integral over $\C
[X_i']$ and $q_i$ is finite. Consider the commutative diagram
\[ \begin{array}{ccc}
X & \stackrel{\rho}{\rightarrow} & X_1\times X_2\\ \, \, \, \,
\downarrow p && \, \,  \, \, \downarrow q \\ X' &
\stackrel{\rho'}{\rightarrow} &
X_1'\times X_2'\\
\end{array} \]
\noindent where $\rho =(\rho_1,\rho_2), \,  \rho'
=(\rho_1',\rho_2')$, and $q=(q_1,q_2)$. Let $Z$ (resp. $Z'$) be
the closure of $\rho (X)$ in $X_1\times X_2$ (resp. $\rho' (X')$
in $X_1'\times X_2'$). By Proposition \ref{app.4} $\rho: X\to Z$
is an isomorphism and, therefore, by Lemma 3.6 in \cite{KaKu1}
$\rho' : X' \to Z'$ is birational finite. Since the statement of
this Lemma is equivalent to the fact that $\rho'$ is an
isomorphism, it suffices to prove $\rho'$ is a holomorphic
embedding.

Consider an orbit $O \subset X$ of $H_1$ and set $O'=p(O)$, $O_2'=\rho_2'
(O')$. Each of these orbits is isomorphic to $\C_+$ and,
therefore, the $H_1$-equivariant finite morphisms $p|_O : O \to
O'$ and $\rho_2'|_{O'}: O' \to O_2'$ must be isomorphisms. Thus
one has a regular function on $X_2'$ whose restriction yields a
coordinate on $O' \simeq O_2' \subset X_2'$. Since locally $X'$ is
biholomorphic to $U \times O'$ where $U$ is an open subset of
$\rho_1' (X')\subset X_1'$ we see that $\rho' : X' \to X_1' \times
X_2'$ is a local holomorphic embedding, i.e. it remains to show
that $\rho'$ is injective.  For any $x \in X$ set $x'=p(x)$, and
$x_i'=\rho_i' (x')$. Assume that $x$ and $y\in X $ are such that
$(x_1',x_2' ) =(y_1',y_2')$. Arguing as in Lemma 3.6\footnote{It
is shown in that Lemma that $\rho_j (p^{-1} (x'))=q_j^{-1} (x_j')$
for a general point $x \in X$ and to adjust the argument to the
present situation one needs it to be true for every point in $X$,
but this follows, of course, from continuity and finiteness of
$q_j$.} in \cite{KaKu1} we can suppose that $x$ and $y$ belong to
the same fiber of $\rho_1$ that is, by assumption, an $H_1$-orbit
$O$. Since $\rho_2'|_{O'}: O' \to O_2'$ is an isomorphism we have
$x'=y'$ which implies the desired conclusion.

\end{proof}

\begin{lemma}\label{app.8} Let the assumption of Lemma
\ref{app.6} hold with one exception: instead of the finiteness of
$p$ we suppose that there are a surjective quasi-finite morphism
$r : S \to S'$ of normal affine algebraic varieties equipped with
trivial $\Gamma$-actions and a surjective $\Gamma$-equivariant
morphism $\varrho' : X' \to S'$ such that $X$ is isomorphic to
fibred product $X' \times_{S'}S$ with $p : X\to X'$ being the
natural projection (i.e. $p$ is surjective quasi-finite). Then the
conclusion of Lemma \ref{app.6} remains valid.

\end{lemma}

\begin{proof}

By construction, $X_i =X_i' \times_{S'} S$. Thus we have the
following commutative diagram
\[ \begin{array}{ccccccc}
X & \stackrel{{ \rho}}{\rightarrow} & ({X}_1'\times
{X}_2')\times_{(S' \times S')} (S \times S) & \stackrel{{( \tau
,\tau ) }}{\rightarrow}
 & S \times S\\
\, \, \, \, \downarrow {p} &&
\, \,  \, \, \downarrow {q} &&  \, \,  \, \, \,  \, \, \, \, \, \downarrow {(r, r)}\\
X' & \stackrel{\rho'}{\rightarrow} &
X_1'\times X_2'  \,  \, & \stackrel{{(\tau' , \tau' )}}{\rightarrow} & \, \, S'\times S'.\\
\end{array} \]
Set $Z$ (resp. $Z'$) equal to the closure of $\rho (X)$ in
$X_1\times X_2$ (resp.  $\rho' (X')$ in $X_1'\times X_2'$) and $D
\simeq S $ (resp. $D' \simeq S'$) be the diagonal subset in
$S\times S$ (resp. $S' \times S')$. 
Since $X= X' \times_{S'} S$ we see that $Z=Z' \times_{D'} D$.

Assume that $\rho' (x')=\rho' (y') =: z'$ for some $x', y' \in X'$. Then by the
commutativity of the diagram we have also  $\varrho' (x')=\varrho' (y') = :s'$.
Since $r$ is surjective $r (s)=s'$ for some $s \in S$. Thus the  elements
$(x', s)$ and $(y',s)$ of   $X' \times_{S'} S$ go to the same element $(z',s')$
of   $Z' \times_{D'} D$ under morphism $\rho$. By
Lemma \ref{app.4} $\rho : X \to Z$ is an isomorphism and therefore $x'=y'$.
Hence $\rho'
: X' \to Z'$ is bijective\footnote{Note that this (slightly modified) argument
provides a much simpler proof of Lemma 3.7 in \cite{KaKu1} where
assuming that $\rho$ is birational finite one needs to show that
$\rho'$ is such.}. It was shown in the proof of Lemma \ref{app.6}
that $\rho'$ is locally biholomorphic, i.e. it is an isomorphism
which implies the desired conclusion.

\end{proof}

\providecommand{\bysame}{\leavevmode\hboxto3em{\hrulefill}\thinspace}

\end{document}